\documentclass{elsarticle}
\usepackage{amsmath}

\usepackage{amssymb}

\journal{Appl.\ Numerical Math.}

\begin{document}

\begin{frontmatter}



\title{Numerical Solution of an Extra-wide Angle Parabolic Equation through Diagonalization of a 1-D Indefinite Schr\"{o}dinger Operator with a Piecewise Constant Potential}


\author[USM]{Sarah D. Wright\corref{cor1}}
\cortext[cor1]{Corresponding author.} 
\author[USM]{James V. Lambers}

\address[USM]{School of Mathematics and Natural Sciences, The University of Southern Mississippi,
            118 College Dr \#5043, 
            Hattiesburg, MS
            39406, USA}

\begin{abstract}
We present a numerical method for computing the solution of a partial differential equation (PDE) for modeling acoustic pressure, known as an extra-wide angle parabolic equation, that features the square root of a differential operator. The differential operator is the negative of an indefinite Schr\"{o}dinger operator with a  piecewise constant potential. This work primarily deals with the 3-piece case; however, a generalization is made the case of an arbitrary number of pieces. Through restriction to a judiciously chosen lower-dimensional subspace, approximate eigenfunctions are used to obtain estimates for the eigenvalues of the operator. Then, the estimated eigenvalues are used as initial guesses for the Secant Method to find the exact eigenvalues, up to roundoff error. An eigenfunction expansion of the solution is then constructed. The computational expense of obtaining each eigenpair is independent of the grid size. The accuracy, efficiency, and scalability of this method is shown through numerical experiments and comparisons with other methods.
\end{abstract}




\end{frontmatter}


\section{Introduction} \label{sec_intro}

The problem to be considered in this paper is an extra-wide angle parabolic equation (termed EWAPE2 in \cite{erdc}), a PDE that is applicable to sound wave propagation where backscattered energy is negligible \cite{gilbert}.
The PDE is expressed in cylindrical coordinates as
\begin{equation}
\frac{du}{dr} = \pm  i\left(\sqrt{\frac{\partial ^2}{\partial z^2} + \alpha^2I}\right)u \hspace{1 cm} 0<z<\pi, \hspace{0.2cm} r>0,
\label{PDE}
\end{equation}
where  $\alpha(r,z)$ is defined to be
\begin{equation}
\alpha=\frac{\omega}{c(r,z)}.
\label{variables}
\end{equation}
Here,  $\omega$ is the angular frequency, and $c(r,z)$ is the sound speed as a function of the range, $r$, and depth, $z$. The range will be treated as a temporal variable for the purpose of numerical solution.

This PDE has many applications in physics \cite{erdc}; our interest is in modeling acoustic pressure in the ocean, in the case where the sound speed is weakly dependent on the range $r$.  In this context, the sound speed has relatively small variation from its average with respect to depth \cite{josie}. Therefore, we treat the sound speed as constant within each of two boundary layers, and within a third interior region. That is, we treat the wavenumber as a piecewise constant coefficient, $\alpha(z)$, over three pieces. Previous work has demonstrated that solution methods using eigenfunction expansion have been effective in solving PDEs with piecewise constant coefficients \cite{Long}.  In this paper, highly accurate eigenvalues and eigenfunctions will be computed, and then used to construct the solution of our version of the EWAPE2, in a way that sidesteps the difficulties encountered by numerical methods due to the
need to approximate the square root of a differential operator.
Although this work is connected to previous work on related eigenvalue problems \cite{Long,elyse}, a key difference is that two distinct approaches are required to obtain reasonably accurate initial guesses 
for the eigenvalues, that can then be improved via iteration--one for low-frequency eigenfunctions, and one for high-frequency.  By exploiting the orthogonality of eigenspaces, an eigenfunction of either kind can be
obtained in $O(1)$ floating-point operations.


The outline of the paper is as follows. In Section 2, we will discuss the background of the problem in further detail. Section 3 will describe the approach to finding eigenfunctions. Section 4 will discuss the solution of the PDE through an eigenfunction expansion.  Section 5 is a compilation of all of the numerical results obtained using our algorithm. This will include evaluation of accuracy and efficiency and comparison of our results to those of other numerical methods. Section 6 will describe a generalization of our algorithm to an arbitrary piecewise constant wavenumber.  Section 7 will conclude the paper with final thoughts regarding the proposed algorithm, and possible future work.

\section{Background} \label{sec_background}
The Helmholtz equation is an eigenvalue problem for the Laplacian operator and has many applications in physics such as seismology, acoustics, and electromagnetic radiation \cite{helmholtz}. 
In cylindrical coordinates, it has the form
$$\left(\frac{\partial^2}{\partial r^2}+\frac{1}{r}\frac{\partial}{\partial r}+\frac{\partial^2}{\partial z^2}+\alpha^2\right)P=0.$$
This equation can be factored into two equations that model either forward propagation or backward propagation. The resulting PDE, also known as the ``one-way equation", is defined as
\begin{equation} \label{eq:thePDE}
\frac{du}{dr} = \pm  i \sqrt{L} u, \quad 0<z<\pi, \quad r>0,
\end{equation}
where the spatial domain (that is, the range of depths) is chosen to be $(0,\pi)$ solely for notational convenience, and the differential operator $L$ is defined by
\begin{equation} \label{eq:theL}
L = \frac{\partial ^2}{\partial z^2} + \alpha^2I.
\end{equation}
This one-way equation only models waves traveling in one direction and so we must assume that any propagation in the opposite direction, such as echoes, is negligible.

We also impose the ``initial" condition
\begin{equation} \label{eq:theIC}
u(z,0)=f(z), \quad 0<z<\pi,
\end{equation}
and homogeneous Dirichlet boundary conditions. 
The solution $u$ to this one-way equation is related to the solution, $P$, of the Helmholtz equation by the change of variable 
$$u=\sqrt{r}P.$$

We are assuming that the sound speed depends only on the depth and not range, to justify the factoring of the Helmholtz equation \cite{gilbert}. 
Because the sound speed flucutates very little relative to its average, we can approach this problem 
as if it were nearly piecewise constant. Therefore, for our simulation, we will model $\alpha^2$ as a piecewise constant function consisting of three pieces, where the first and last piece come from the wavenumber at the bounding depths. The average of the wavenumber over the interior depths will be used for the middle piece. 

The one-way equation (\ref{eq:thePDE}) includes the square root of a negated 1-D Schr\"{o}dinger operator. A general Schr\"{o}dinger operator is the sum of the kinetic energies, represented as a negated Laplacian, and the potential energies, accounted for by $V$ \cite{Schrodinger}:
$$H=-\left(\sum_{j=1}^{m}{\frac{\partial^2}{\partial x_j^2}}\right)+V.$$
Our operator $L$ in the one-way equation is a scalar multiple of a Schr\"{o}dinger operator with a negative potential; such Schr\"{o}dinger operators were studied in \cite{fefferman}.
Taking the square root of $L$ poses substantial difficulties for numerical methods.  For a narrow-angle parabolic equation,
$\sqrt{L}$ is approximated by using a first-order Taylor expansion of 
$$\alpha_0\sqrt{1+\ell}$$
where
$$\ell=\frac{L}{\alpha_0^2}-1.$$
This was originally proposed by Claerbout for extrapolation of seismic waves \cite{claerbout}. In order to obtain a sufficiently accurate approximation, Taylor expansion is used, which yields
$$\sqrt{1+\ell} \approx 1+\frac{\ell}{2}.$$
By limiting the Taylor expansion to only the first two terms, this approach sacrifices accuracy for efficiency. This approach will not work well for our problem because the matrix $A$ obtained by discretizing $L$ is ill-conditioned.  For a wide-angle parabolic equation, a rational approximation is used instead \cite{gilbert},
but the applicability of this approach is also limited, again due to the ill-conditioning of this matrix.

Another approach is to use the Schur decomposition $A=QDQ^H$, where 
$Q$ is unitary, and $D$ is diagonal, as $A$ is symmetric. It follows that $\sqrt{A}=Q\sqrt{D}Q^H$. 
This method is very expensive, as it requires approximately $28\frac{1}{3}n^3$ flops \cite{higham}, 
where $A$ is $n\times n$. A far more efficient approach would be to compute accurate eigenvalues and eigenfunctions of the operator directly, without relying on any spatial discretization.  We will develop such an approach in this paper.

\subsection{Previous Work}
Similar to previous work concerning eigenfunction expansion \cite{Long,elyse}, the derivation of the proposed algorithm will begin by defining the eigenfunctions as combinations of sine and cosine functions. Each eigenfunction must satisfy boundary conditions, and continuity requirements at the interfaces between pieces. The system of equations arising from these conditions can then be reduced by eliminating most parameters. Once the system is reduced to a single nonlinear equation, we can solve for the remaining parameter using the Secant Method \cite{burden}. This remaining variable is the frequency corresponding to the third piece of the eigenfunction. 
In order to use The Secant Method, we will seek simple approximations of each eigenvalue, to serve 
as initial guesses, as was done in \cite{Long,elyse}, but as the differential operator is of a different form,
a substantially different approach to obtaining these initial guesses will be required.

\subsection{Behavior of Eigenfunctions}
To guide the development of our algorithm, we use {\sc Matlab} to compute eigenvalues and eigenvectors of a matrix
that discretizes our operator $L$, so that we may understand their behavior.
\begin{figure}[ht]
\begin{center}
\includegraphics[width=2.5in]{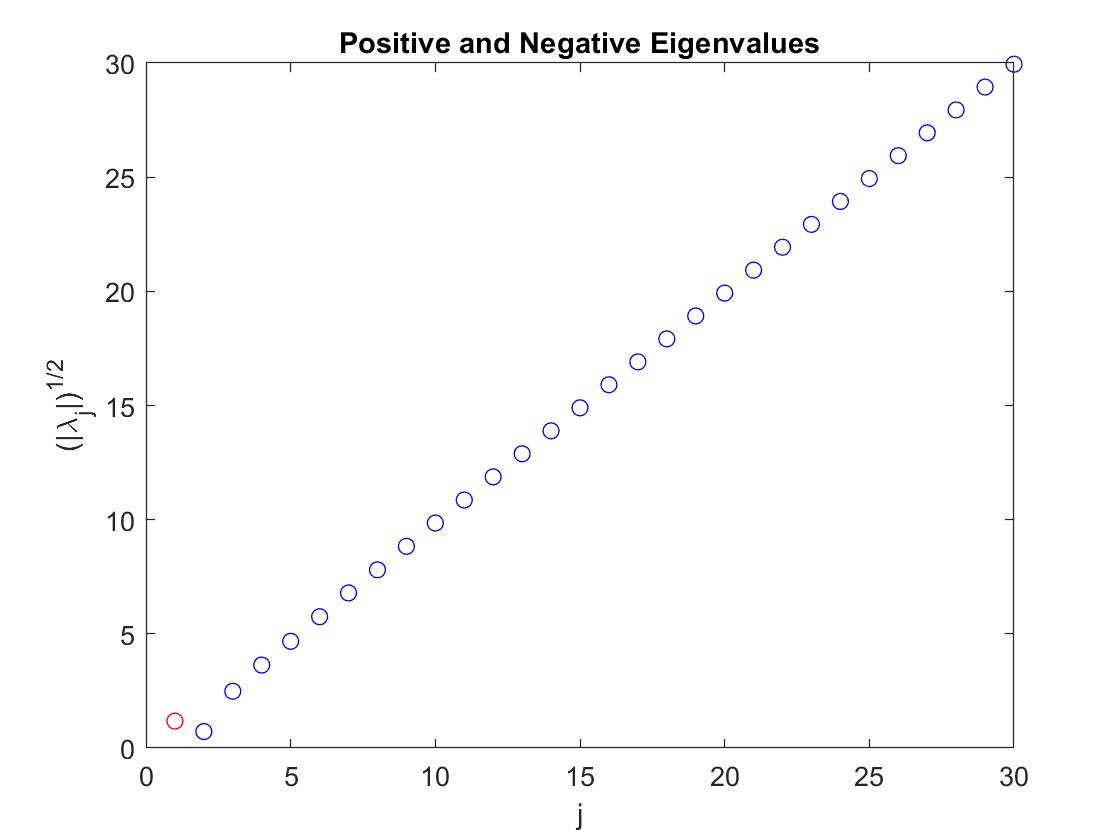}
\end{center}
\caption{Absolute values of all eigenvalues of ``one way" equation with $\alpha_1=[2, 1, 2]$, $\rho=[1/3, 2/3]$, and $N=2048$. The positive eigenvalues are seen in red and the eigenvalues that correspond to negative eigenvalues are blue.}
\label{A1}
\end{figure}
\begin{figure}[ht]
\begin{center}
\includegraphics[width=2.5in]{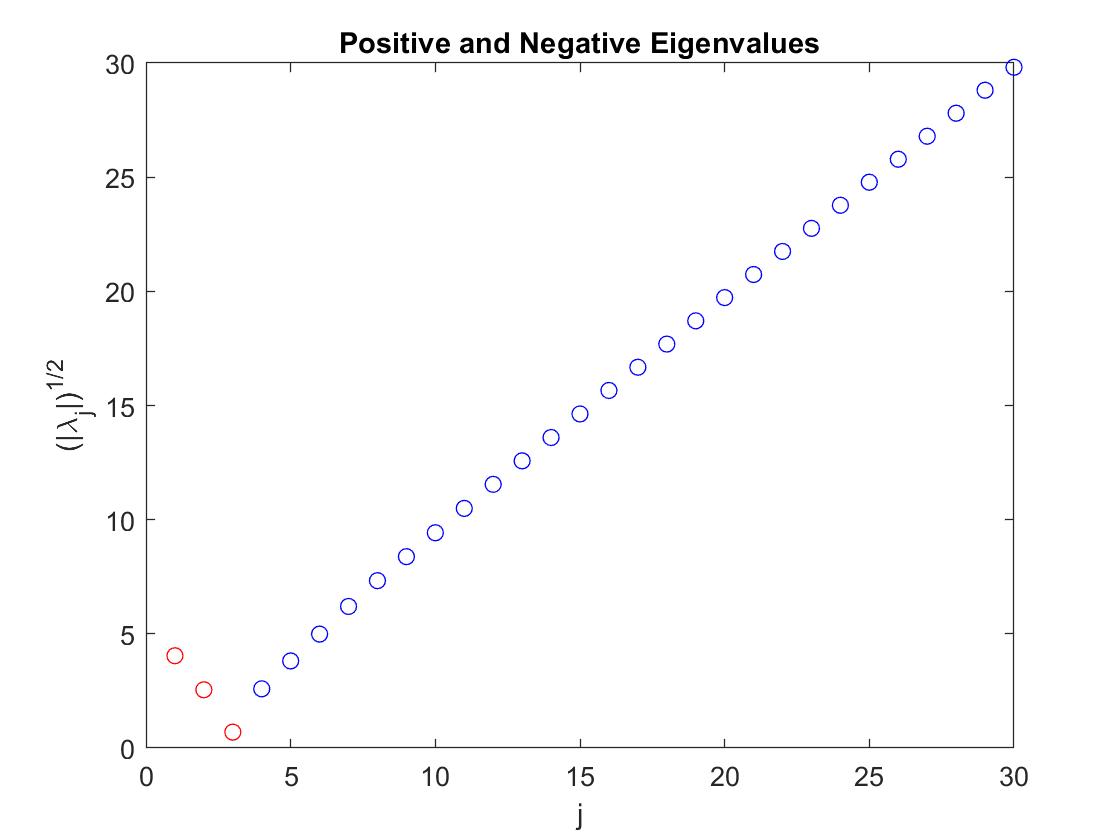}
\end{center}
\caption{Absolute values of all eigenvalues of ``one way" equation with $\alpha_2=[1, 3, 5]$, $\rho=[1/4, 3/4]$, and $N=2048$. The positive eigenvalues are seen in red and the eigenvalues that correspond to negative eigenvalues are blue.}
\label{A2}
\end{figure}
\begin{figure}[ht]
\begin{center}
\includegraphics[width=2.5in]{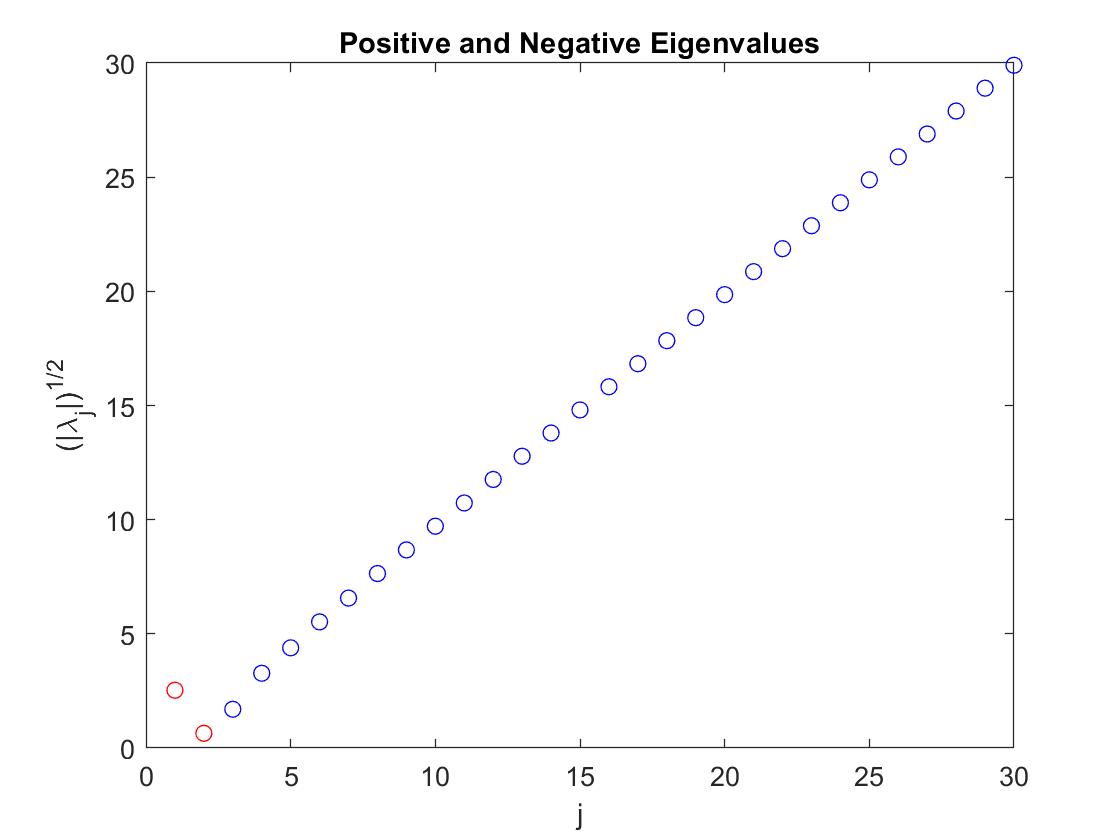}
\end{center}
\caption{Absolute values of all eigenvalues of ``one way" equation with $\alpha_3=[2, 1, 3]$, $\rho=[1/4, 1/2]$, and $N=2048$. The positive eigenvalues are seen in red and the eigenvalues that correspond to negative eigenvalues are blue.}
\label{A3}
\end{figure}
Figures \ref{A1}, \ref{A2}, and \ref{A3} show the graphs of the absolute value of the smallest several eigenvalues. The first few are positive eigenvalues that are bounded in magnitude; however, the remaining eigenvalues are negative and are unbounded. 

\begin{figure}[ht]
\begin{center}
\includegraphics[width=2in]{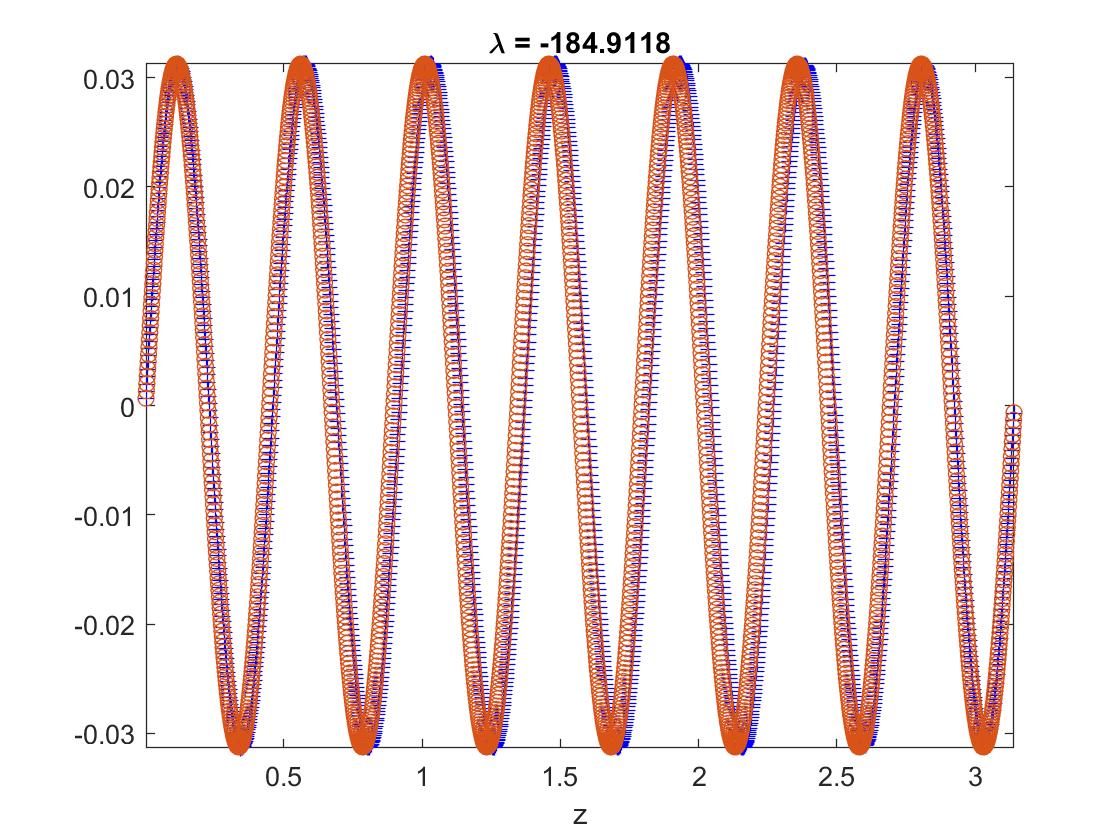} 
\includegraphics[width=2in]{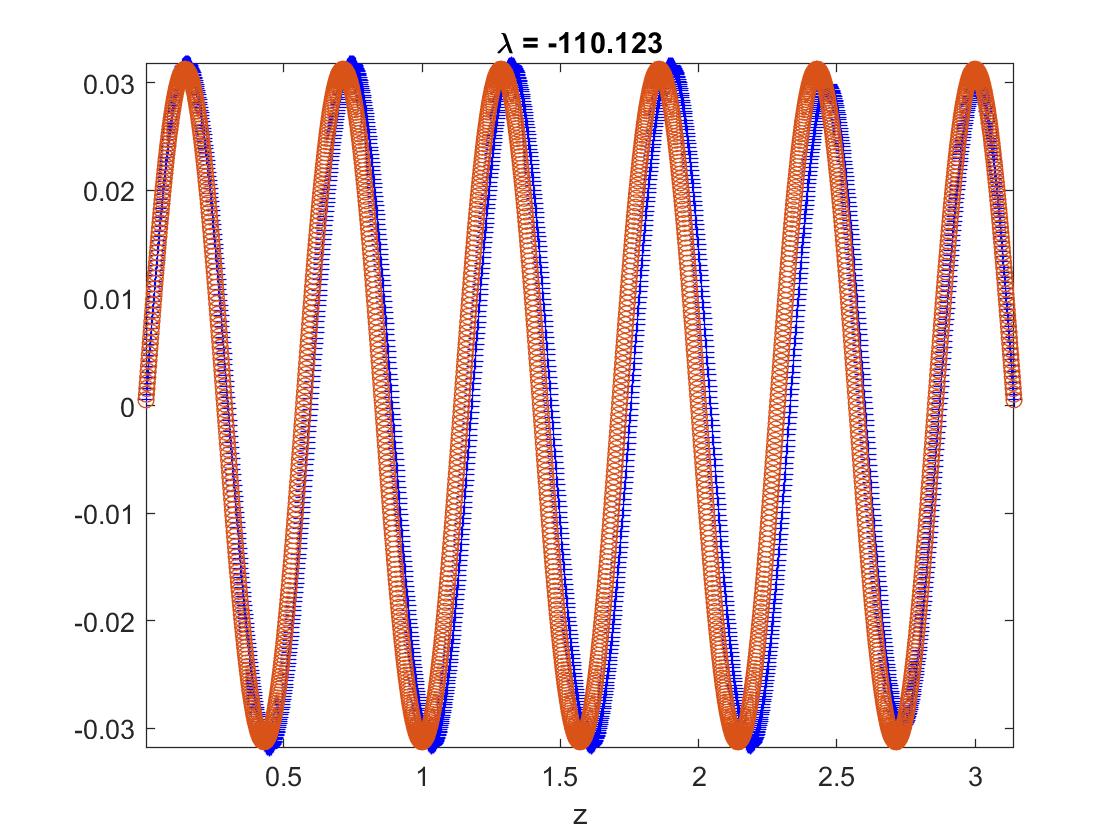}\\
\includegraphics[width=2in]{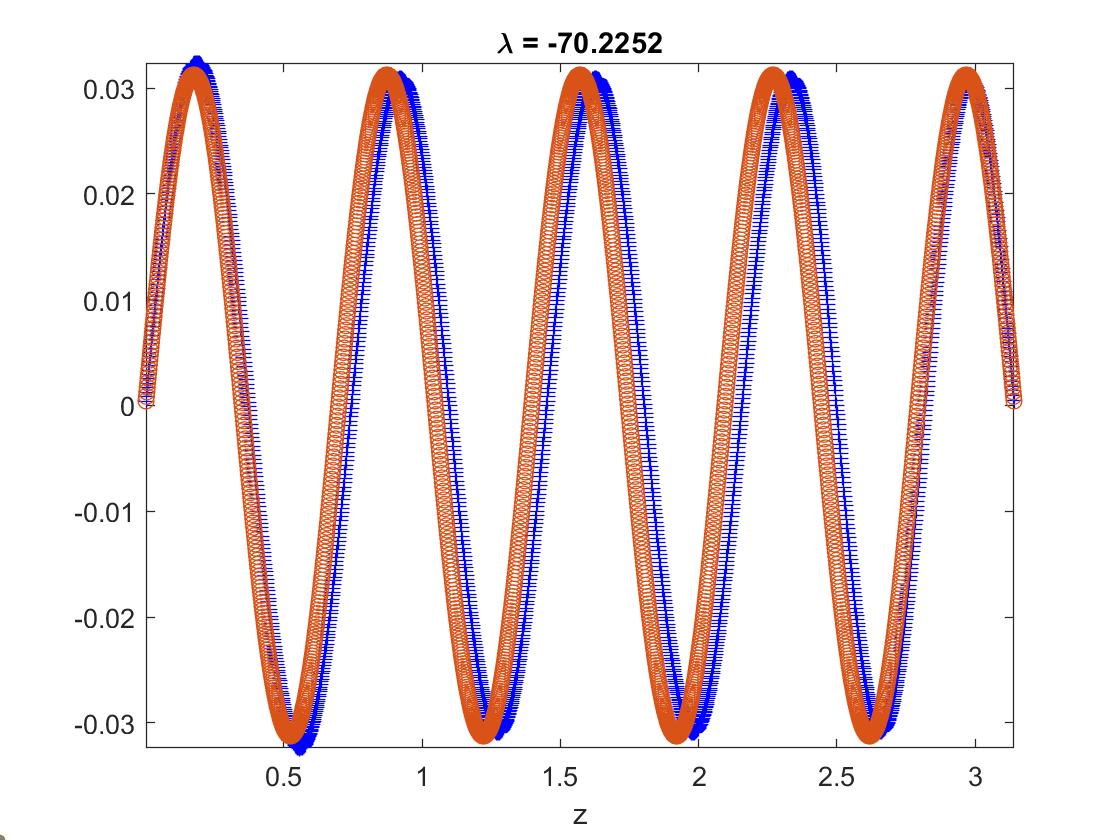}
\includegraphics[width=2in]{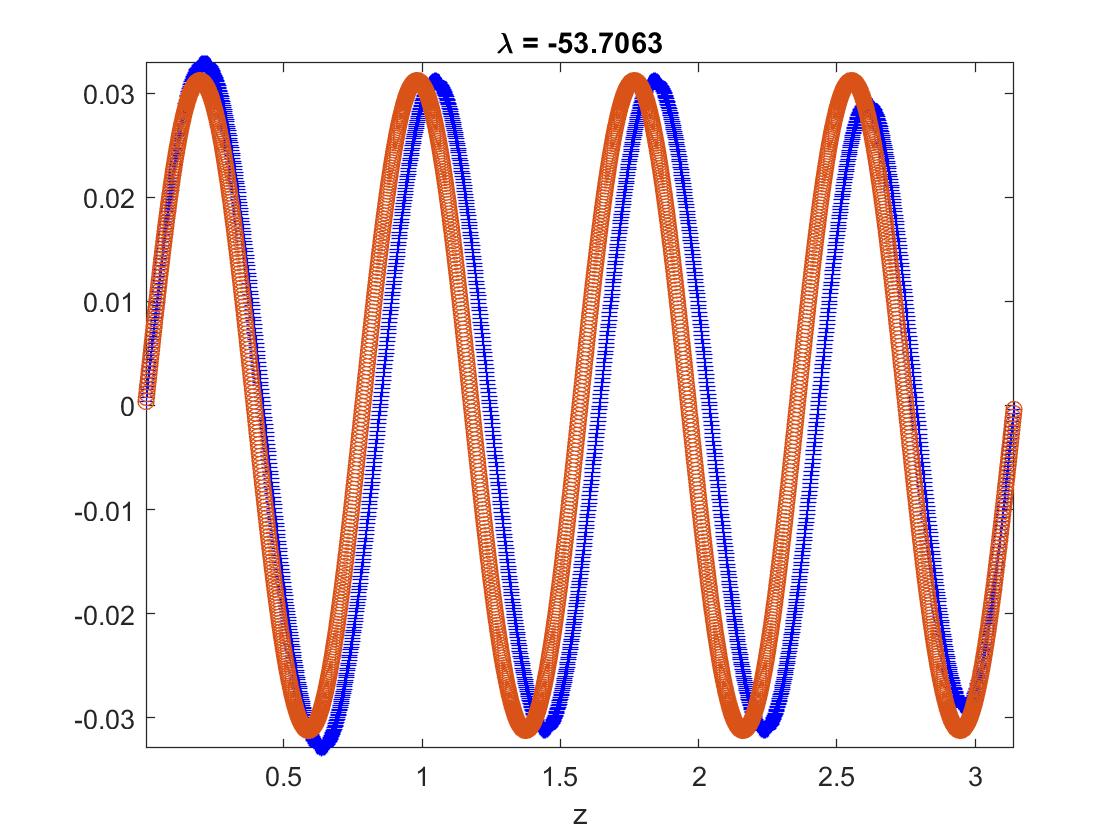}
\end{center}
\caption{Eigenfunctions that correspond to negative  eigenvalues compared to the function $\sin(jz)$.}
\label{Negative}
\end{figure}

\begin{figure}[ht]
\begin{center}
\includegraphics[width=2in]{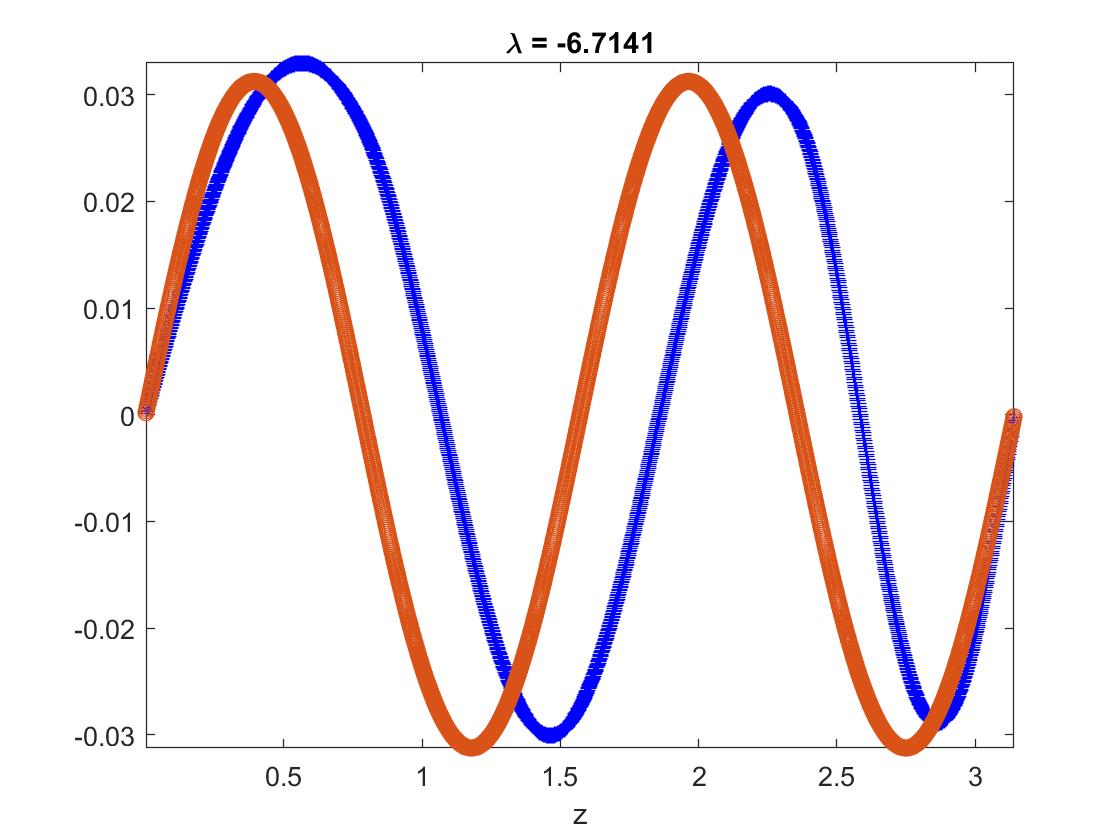}
\includegraphics[width=2in]{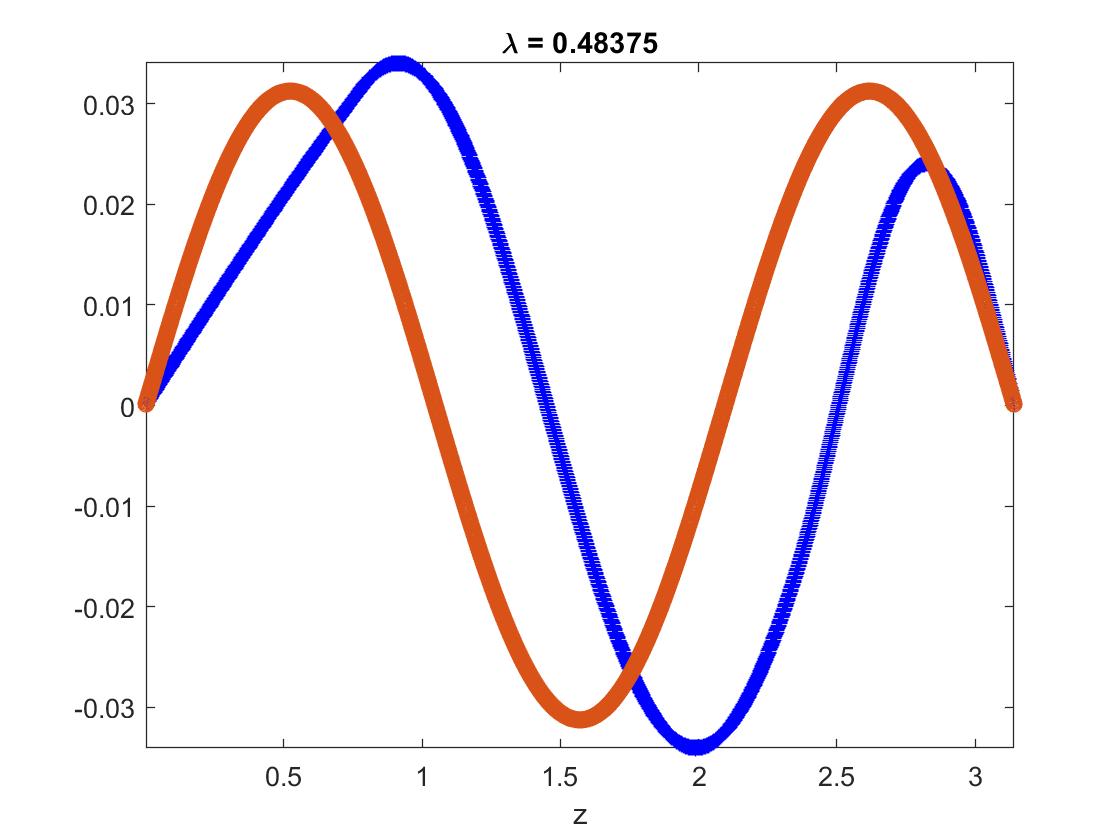}\\
\includegraphics[width=2in]{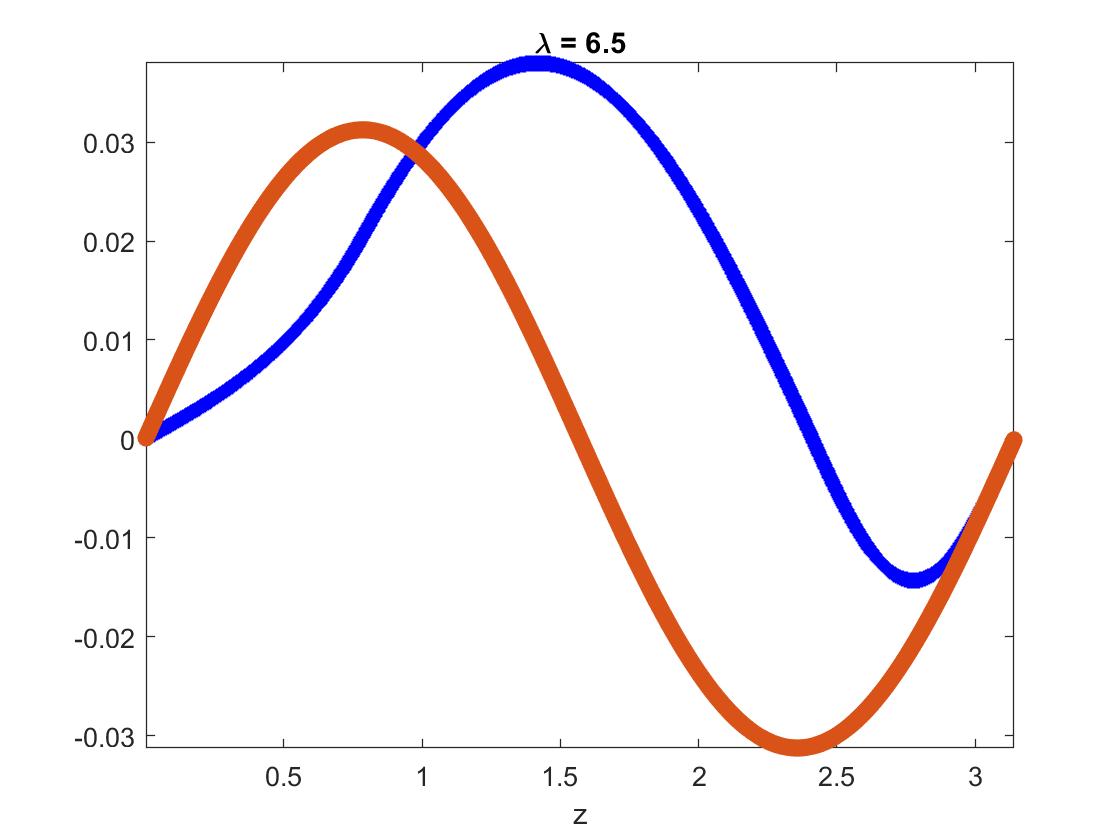}
\includegraphics[width=2in]{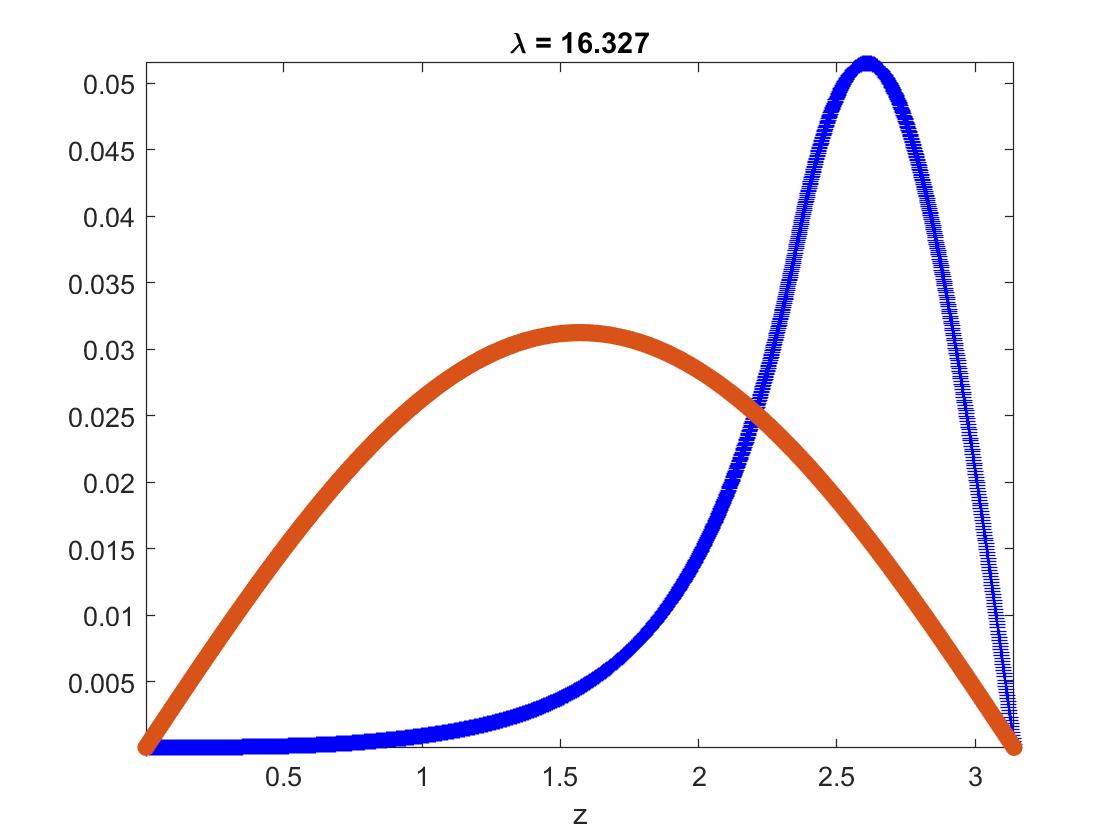}
\end{center}
\caption{Eigenfunctions that correspond to positive eigenvalues or are close to 0 compared to the function $\sin(jz)$.}
\label{positive}
\end{figure}
Each plot within Figure \ref{Negative} displays an eigenfunction in blue and the function $\sin(jz)$, for some $j$, in red. We can see that eigenfunctions that correspond to larger eigenvalues can be approximated by a sine function; however, as the eigenvalues approach 0, they become more difficult to approximate using such a sine function. In fact, when the eigenvalues become positive, they nowhere resemble a sine function, as we can see from Figure \ref{positive}, so we must therefore develop a different approach for approximating these eigenfunctions.

\section{Methodology}

The eigenfunctions for the piecewise constant coefficient case, with three pieces, will have the form
$$V_j(z)=
\left \{
  \begin{array}{ll}
  V_{j1}(z)=A_{j1}\cos(\omega_{j1}z)+B_{j1}\sin(\omega_{j1}z)	&	0 \leq z<\pi\rho_1, \\
 V_{j2}(z)=A_{j2}\cos(\omega_{j2}z)+B_{j2}\sin(\omega_{j2}z)  &	\pi\rho_1 \leq z<\pi\rho_2,\\
  V_{j3}(z)=A_{j3}\cos(\omega_{j3}z)+B_{j3}\sin(\omega_{j3}z)	&	\pi\rho_2 \leq z<\pi,
    \end{array}
\right.$$
where $\omega_{j1}, \omega_{j2},\omega_{j3}$ are the unknown frequencies and $A_{j1}, A_{j2}, A_{j3}$ and $B_{j1},B_{j2}, B_{j3}$ are related to the unknown amplitudes and phase shifts.
A system of linear equations that characterizes the eigenfunction comes from imposing conditions on $V_{j}(z)$. These conditions are homogeneous Dirichlet boundary conditions and continuity conditions. The continuity conditions require that the one-sided limits of eigenfunctions at the interfaces are equal, as well as those of their first derivatives \cite{gottlieb}. We will also set conditions for how the frequencies $\omega_{j1}$, $\omega_{j2}$, and $\omega_{j3}$ relate to each other. 
\begin{itemize}
\item{Dirichlet Boundary Conditions:} 
\begin{equation}  V_{j1}(0)=0, \label{BC1} \end{equation}
\begin{equation}  V_{jn}(\pi)=0. \label{BC2} \end{equation}

\item{Continuity:} 
\begin{equation}  V_{j1}(\pi\rho_1)=V_{j2}(\pi\rho_1), \label{C1} \end{equation}
\begin{equation} V'_{j1}(\pi\rho_1)=V'_{j2}(\pi\rho_1), \label{C2} \end{equation}
\begin{equation}  V_{j2}(\pi\rho_2)=V_{j3}(\pi\rho_2), \label{C3} \end{equation}
\begin{equation} V'_{j2}(\pi\rho_2)=V'_{j3}(\pi\rho_2). \label{C4} \end{equation}
\end{itemize}
This leads to a homogeneous system of equations for the coefficients $A_i$ and $B_i$. 
From the boundary conditions we have
\begin{eqnarray} 
A_{j1}\cos(\omega_{j1}(0))+B_{j1}\sin(\omega_{j1}(0)) \,=\, A_{j1} & = & 0,\\
A_{j3}\cos(\omega_{j3}(\pi))+B_{j3}\sin(\omega_{j3}(\pi))  &=&  0.
\end{eqnarray}
From the continuity requirements, we have
\begin{eqnarray}
A_{j1}\cos(\omega_{j1}\pi\rho_1) + B_{j1}\sin(\omega_{j1}\pi\rho_1) & = & A_{j2}\cos(\omega_{j2}\pi\rho_1) \nonumber\\
&&+B_{j2}\sin(\omega_{j2}\pi\rho_1) \\
-\omega_{j1}A_{j1}\sin(\omega_{j1}\pi\rho_1) + \omega_{j1}B_{j1}\cos(\omega_{j1}\pi\rho_1) & = & -\omega_{j2}A_{j2}\sin(\omega_{j2}\pi\rho_1) \nonumber \\
&& +\omega_{j2}B_{j2}\cos(\omega_{j2}\pi\rho_1) \\
A_{j2}\cos(\omega_{j2}\pi\rho_2) + B_{j2}\sin(\omega_{j2}\pi\rho_2) & = & A_{j3}\cos(\omega_{j3}\pi\rho_2) \nonumber\\
&& + B_{j3}\sin(\omega_{j3}\pi\rho_2) \\
-\omega_{j2}A_{j2}\sin(\omega_{j2}\pi\rho_2) + \omega_{j2}B_{j2}\cos(\omega_{j2}\pi\rho_2) & = & -\omega_{j3}A_{j3}\sin(\omega_{j3}\pi\rho_2) \nonumber \\
&& +\omega_{j3}B_{j3}\cos(\omega_{j3}\pi\rho_2). 
\end{eqnarray}

The system of equations from the boundary and continuity requirements results in a homogeneous system of the form $M_j{\bf x}_j={\bf 0},$ where 
\begin{align*}
{\bf x}=\begin{bmatrix}
x_1\\
{\bf x}_2\\
\end{bmatrix}, \quad  x_{j1}=B_{j1}, \quad {\bf x}_{j2}=\begin{bmatrix}
A_{j2}\\
B_{j2}\\
A_{j3}\\
B_{j3}\\
\end{bmatrix},
\end{align*}
and $B_{j1}, A_{j2}, B_{j2}, A_{j3}, B_{j3}$ are the unknowns of the system.
The matrix $M$ has the $2\times 2$ block structure
$$M_j= \left[
\begin{array}{cc}
0 & B_j \\
C_j & D_j \\
\end{array}
\right]$$
where
$$
B_j= \left[
\begin{array}{cccc}
0 & 0 & \cos(\pi\omega_{j3}) & \sin(\pi\omega_{j3})
\end{array}
\right], \quad
C_j= \left[
\begin{array}{c}
\sin(\pi\rho_1\omega_{j1}) \\
 \omega_{j1}\cos(\pi\rho_1\omega_{j1}) \\
 0\\ 
0 \\
\end{array}
\right],$$
$$D_j= \left[
\begin{array}{cccc}
-\cos(\pi\rho_1\omega_{j2}) & -\sin(\pi\rho_1\omega_{j2}) & 0 & 0 \\
\omega_{j2}(\sin(\pi\rho_1\omega_{j2})) &  -\omega_{j2}
(\cos(\pi\rho_1\omega_{j2})) &  0 & 0 \\
\cos(\pi\rho_2\omega_{j2}) & \sin(\pi\rho_2\omega_{j2}) & -\cos(\omega_{j3}\pi\rho_2) & -\sin(\omega_{j3}\pi\rho_2) \\
-\omega_{j2}(\sin(\pi\rho_2\omega_{j2})) & \omega_{j2}(\cos(\pi\rho_2\omega_{j2}) & \omega_{j3}\sin(\omega_{j3}\pi\rho_2) & -\omega_{j3}\cos(\omega_{j3}\pi\rho_2) \\
\end{array}
\right].$$
The entries of this matrix come from the coefficients of $A_{ji}$ and $B_{ji}$ from the system of Eqns. (\ref{BC1})-(\ref{C4}).

On each piece, $V_{j}(z)$ is an eigenfunction of a constant-coefficient operator. For $V_j(z)$ to be an eigenfunction on the entire spatial domain $[0,\pi]$, the eigenvalues $\lambda_{j1}, \lambda_{j2}, \lambda_{j3}$ corresponding the three pieces must be equal. This requirement, as well as the equation $\lambda_{ji}=-\omega_i^2+\alpha_i^2$, yields the relationship between the $\omega_{ji}$ values, for $i=1,2,3$:
$$\lambda_{j1}=\lambda_{j3}$$
$$-\omega_{j1}^2+\alpha_1^2=-\omega_{j3}^2+\alpha_3^2$$
which can be rewritten as 
\begin{eqnarray}
\omega_{j1} & = & \sqrt{\omega_{j3}^2 - \alpha_3^2 +\alpha_1^2}
\end{eqnarray}
and 
$$\lambda_{j2}=\lambda_{j3}$$
$$-\omega_{j2}^2+\alpha_2^2=-\omega_{j3}^2+\alpha_3^2$$
which can be rewritten as
\begin{eqnarray}
\omega_{j2} & = & \sqrt{\omega_{j3}^2 - \alpha_3^2 + \alpha_2^2}
\end{eqnarray}
so that $\omega_{j1}$ and $\omega_{j2}$ are eliminated.

To obtain an eigenfunction, we must find a value for $\omega_{j3}$ such that the matrix $M_j$ is singular.  We have the following equation for the determinant, 
\begin{equation} 
\det(M_j)=\det(-B_j D_j^{-1}C_j)\det(D_j).
\end{equation}
In an effort to simplify the process of solving for $\omega_{j3}$, we can prove that the determinant of the block $D_j$ is always non-zero; in fact, $\det(D_j)=\omega_{j3}\omega_{j2}$. Therefore we will proceed to find a value of $\omega_{j3}$ by solving  $\det(-B_jD_j^{-1}C_j)=0$. 
The Secant Method will be used to determine the value for $\omega_{j3}$; however, in order to use the Secant Method, we must have two initial guesses.

For sufficiently large $j$ values, we have seen that $\sin(jz)$ is a good approximation of an eigenfunction of $L$. In order to obtain an approximate $\omega_{j3}$ value, we will use the Rayleigh Quotient
\begin{equation}
R(L,F_j)=\frac{\langle F_j, LF_j \rangle}{\langle F_j, F_j \rangle}
\end{equation}
where we are using the standard inner product of real-valued functions on $[0,\pi]$,
\begin{equation}
\langle f, g\rangle=\int_{0}^{\pi}f(z)g(z)\,dz.
\end{equation}
We will see that the Rayleigh Quotient is an accurate approximation of the exact eigenvalue $\lambda$ for sufficiently large $j$. This yields the first initial guess of the eigenvalue,
\begin{equation}
\lambda_j^{(0)}=\frac{\langle \sin(j\cdot),L\sin(j\cdot) \rangle}{\langle \sin(j\cdot), \sin(j\cdot) \rangle}.
\end{equation}
Because the exact eigenvalue $\lambda_j$ has the form
\begin{equation}
\lambda_j=-\omega_{ji}^2+\alpha_i^2, \label{Eig} \quad i=1,2,3,
\end{equation}
we can readily obtain a first initial guess for $\omega_{j3}$. We can set the right side of this equation equal to the Rayleigh Quotient, as follows:
\begin{equation}
\frac{\langle \sin(j\cdot),L\sin(j\cdot) \rangle}{\langle \sin(j\cdot), \sin(j\cdot) \rangle}= -\omega_i^2+\alpha_i^2, \label{IG}
\end{equation}
and solve for $\omega_{ji},$ which will be our first initial guess to use for the Secant Method.

To find an initial guess for $\omega_{j3}$ for smaller values of the index $j$, we will need another approach. As Figure \ref{positive} illustrated, the function $\sin(jz)$ is not a good approximation of our eigenfunction for small $j$. For $j=1,2,\ldots,J$, we define $\phi_j(z)=\sin(jz)$.
Here, $J$ represents a hypothetical ``cutoff value", or the largest value of $j$ for which $\sin(jz)$ is not a
sufficiently accurate approximate eigenfunction.  
We then define the subspace spanned by these inaccurate approximations,
$$\Phi_{J}=\textrm{span}\{\phi_1, \dots, \phi_J\}.$$
The orthogonal projection $P_J$ of a function $u$ onto $\Phi_J$ is defined by
$$P_Ju=\sum_{j=1}^{J}\frac{\phi_j\langle \phi_j, u\rangle}{\langle \phi_j, \phi_j \rangle}.$$
Then, the restriction of the operator $L$ to \textbf{$\phi_J$} is defined by
$$L_J=P_JLP_J.$$
The matrix $L_0$ of $L_J$ has entries
\begin{eqnarray*}
\left[ L_0 \right]_{ij}&=&\frac{\langle\phi_i, L\phi_j \rangle}{\langle \phi_i, \phi_i \rangle}\\
&=&\frac{\int_{0}^{\pi}\phi_i(z)\left(\frac{\partial^2}{\partial(z)^2}+\alpha^2(z)\right)\phi_j(z)\,dz}{\sqrt{\int_{0}^{\pi}\phi_i^2(z)\,dz\int_{0}^{\pi}\phi_j^2(z)\,dz}}, \quad i,j=1,\ldots,J,
\end{eqnarray*}

First, we will construct the matrix entries for $i\neq j$.
From $\phi_n(z)=\sin(nz)$ and $L\phi_n=-n^2\sin(nz)+\alpha(z)^2\sin(nz)$, we obtain
\begin{eqnarray*}
\left[L_0 \right]_{ij}&=&\frac{\int_{0}^{\pi}\left(\sin(iz)\right)\left(-j^2\sin(jz)+\alpha(z)^2\sin(jz)\right)dz}{\int_{0}^{\pi}\sin(jz)\sin(jz)\,dz}\\
&=&\frac{2X_{ij}}{\pi}
\end{eqnarray*}
where
\begin{eqnarray}
X_{ij}&=&\frac{\alpha_1^2}{2}\left[ \frac{(i+j)\sin(\rho_1(i\pi-j\pi))-(i-j)\sin(\rho_1(i\pi+j\pi))}{i^2-j^2}\right]+\nonumber\\
&&\frac{\alpha_2^2}{2}\left[ \frac{(i+j)\sin(\rho_2(i\pi-j\pi))-(i+j)\sin(\rho_1(i\pi-j\pi))}{i^2-j^2}\right.\nonumber\\
&&\left.\frac{-(i-j)\sin(\rho_2(i\pi+j\pi))+(i-j)\sin(\rho_1(i\pi+j\pi))}{i^2-j^2}\right]+\nonumber\\
&&\frac{\alpha_3^2}{2}\left[ \frac{-(i+j)\sin(\rho_2(i\pi-j\pi))+(i-j)\sin(\rho_2(i\pi+j\pi))}{i^2-j^2}\right].
\label{Eqn3}
\end{eqnarray}

For the $i=j$ case,
$$\left[L_0 \right]_{ii}=-j^2+\frac{2}{\pi}Y_j,$$
where
\begin{eqnarray}
Y_j&=& \alpha_{1}^2\left[-\frac{\sin(2\pi\rho_{1}j)-2\pi\rho_{1}j}{4j}\right] +\nonumber\\
&&\alpha_{2}^2\left[-\frac{\sin(2\pi\rho_{2}j)-2\pi\rho_{2}j-\sin(2\pi\rho_{1}j)+2\pi\rho_{1}j}{4j}\right] \nonumber\\
 &+& \alpha_{3}^2\left[\frac{\pi}{2} + \frac{\sin(2\pi\rho_{2}j)-2\pi\rho_{2}j}{4j}\right].
\end{eqnarray}

Now that we have formulas to find initial guesses for the values of $j$ for which eigenfunctions can be accurately approximated using $\sin(jz)$, and the matrix $L_0$ corresponding for all other values of $j$, we also need a way to evaluate the cutoff value $J$ that indicates which values of $j$ fall into each category. 
We will test for the values of $j$ for which eigenfunctions cannot be accurately approximated using $\sin(jz)$ using the error formula
$$E_{j}=\frac{\|L\phi_j - \mu_j\phi_j\|^2}{\|L\phi_j\|^2}=\frac{\langle{L\phi_{j},L\phi_{j}}\rangle -\mu_{j}^2 \langle{ \phi_{j},\phi_{j} }\rangle}{\langle{L\phi_{j},L\phi_{j} }\rangle},$$
where  $\mu_{j}=\frac{\langle{ \phi_{j},L\phi_{j} }\rangle}{\langle{\phi_{j},\phi_{j} }\rangle}$ and $\phi_j(z)=\sin(jz)$.
The error formula, once expanded is, 
\begin{eqnarray}
E_j&=& \frac{\displaystyle{ A_j-\frac{2}{\pi}B_j^2}}{\displaystyle{\left(\frac{j^4\pi}{2}\right)-2j^2B_j+A_j}},\label{EF}
\end{eqnarray}
where
%
\begin{eqnarray}
A_j&=&\alpha_{1}^4\left[-\frac{\sin(2\pi\rho_{1}j)-2\pi\rho_{1}j}{4j}\right]\nonumber\\
&& +\alpha_{2}^4\left[-\frac{\sin(2\pi\rho_{2}j)-2\pi\rho_{2}j-\sin(2\pi\rho_{1}j)+2\pi\rho_{1}j}{4j}\right]\nonumber \\
&& + \alpha_{3}^4\left[\frac{\pi}{2} + \frac{\sin(2\pi\rho_{2}j)-2\pi\rho_{2}j}{4j}\right],  \\
B_j&=& \alpha_{1}^2\left[-\frac{\sin(2\pi\rho_{1}j)-2\pi\rho_{1}j}{4j}\right]\nonumber\\
&& +\alpha_{2}^2\left[-\frac{\sin(2\pi\rho_{2}j)-2\pi\rho_{2}j-\sin(2\pi\rho_{1}j)+2\pi\rho_{1}j}{4j}\right]\nonumber \\
&& + \alpha_{3}^2\left[\frac{\pi}{2} + \frac{\sin(2\pi\rho_{2}j)-2\pi\rho_{2}j}{4j}\right].
\end{eqnarray}
This error formula allows us to find the values of $j$ for which $\sin(jz)$ is not a good approximation of the eigenfunction $V_j(z)$, and thus we need to use the matrix $L_0$ to obtain suitable initial guesses. The range of $j$-values that result in a large error from the error formula lead us to a cutoff $j$-value which we call $J$. This is the largest $j$-value that results in a poor approximation using $\sin(jz)$. The matrix $L_0$ is $J \times J$, and its eigenvalues are used to obtain the first initial guess, through the same equation (\ref{Eig}) for each eigenvalue, $\lambda_j$, for $j=1,\ldots , J$.

It is worthwhile to explain why it is sufficient to restrict our attention to the subspace $\Phi_J$ to estimate
the eigenvalues $\lambda_1, \ldots, \lambda_J$.  The exact eigenfunctions $\{ V_j(z) \}_{j=1}^\infty$ are orthogonal, as are
the functions $\sin(jz)$, for $j=1,2,\ldots$,.  As $\sin(jz)$ is an accurate approximation of the eigenfunction corresponding to $\lambda_j$ for $j > J$, it follows that the exact eigenfunctions $V_1(z), \ldots, V_J(z)$ are nearly orthogonal to $\mathrm{span} \{ \sin((J+1)z), \sin((J+2)z), \ldots \}$.  It follows that
a reasonably accurate approximation of $V_j(z)$ can be obtained from the subspace $\Phi_J$, and thus a viable initial guess from the Rayleigh quotient of such an approximation.

 Once the first initial guess is found for each $j$-value, either through the Rayleigh Quotient with $\sin(jz)$ or the eigenvalues of $L_0$, a small perturbation of $\omega_{j3}$ is used for the second initial guess. The initial guesses of $\omega_{j3}$ found for all values of $j$ are passed as input to the Secant Method in order to obtain values of $\omega_{j3}$ that are exact, up to rounding error. Then, each eigenvalues $\lambda_j$ can be obtained by
$$\lambda_j=-\omega_{j3}^2+\alpha_3^2.$$

Then, the matrix $M_j$ is formed, using the computed value of $\omega_{j3}$ such that $M_j$ is singular.
The null space of $M_j$ is found in order to obtain the coefficients of our eigenfunction. 
The equation
$$ \left[
\begin{array}{cc}
 0 & B_j\\
C_j & D_j \\
\end{array}
\right]
\left[
\begin{array}{c}
x_{j1}\\
{\bf x}_{j2}\\
\end{array}
\right]=
\left[
\begin{array}{c}
{0}\\
{\bf 0}\\
\end{array}
\right]$$
indicates how the null space can be obtained, where $B_j$, $C_j$, and $D_j$ are the blocks of the matrix $M_j$. Rewriting this matrix equation, we have
\begin{eqnarray}
B_j{\bf x}_{j2}&=&0, \label{MV1}\\
C_jx_{j1}+D_j{\bf x}_{j2}&=&{\bf 0} \label{MV2}.
\end{eqnarray}
Rearranging equation (\ref{MV2}) yields
\begin{equation}
{\bf x}_{j2}=-D_j^{-1}C_jx_{j1}. \label{MV3}
\end{equation}
From this equation, we see that $x_{j1}$ must be nonzero to ensure that the entire null space vector is nonzero. We will set $x_{j1}=1$, which yields
\begin{equation}
{\bf x}_{j2}=-D_j^{-1}C_j.
\end{equation} 
It can be verified that these choices of $x_{j1}$ and ${\bf x}_2$ satisfy Eqn. (\ref{MV1}), due to $M_j$ being singular. Therefore, we can define the null space as the span of ${\bf x}_j$, where
$${\bf x}_j=\left[\begin{array}{c}
B_{j1}\\
A_{j2}\\
B_{j2}\\
A_{j3}\\
B_{j3}\\
\end{array}\right]=
\left[
\begin{array}{c}
x_{j1}\\
{\bf x}_{j2}\\
\end{array}\right]=
\left[
\begin{array}{c}
1\\
-D_j^{-1}{C}_j\\
\end{array}\right].$$


The last two elements of $C_j$ are $0$; therefore, we only need the first two columns of $D_j^{-1}$. We will use the notation $C_{ik}=\cos(\omega_{ji}\rho_k\pi)$ and $S_{ik}=\sin(\omega_{ji}\rho_k\pi)$, for $i,k=1,2,3$, to simplify the entry formulas. We then have
\begin{eqnarray}
{\bf x}_{j2}&=&-\begin{bmatrix}
-C_{21} & \frac{1}{\omega_{j2}}S_{21}\\
-S_{21} & \frac{1}{\omega_{j2}}C_{21}\\
E & F
\end{bmatrix}
\begin{bmatrix}
S_{11}\\
\omega_{j1}C_{11}\end{bmatrix}
\end{eqnarray}
where $E$ and $F$ are defined by
\begin{eqnarray}
E&=&\begin{bmatrix}
-C_{21}\left[C_{32}C_{22}+\frac{\omega_{j2}}{\omega_{j3}}S_{32}S_{22}\right]-S_{21}\left[C_{32}S_{22}-\frac{\omega_{j2}}{\omega_{j3}}S_{32}C_{22}\right] \\
-C_{21}\left[S_{32}C_{22}-\frac{\omega_{j2}}{\omega_{j3}}C_{32}S_{22}\right]-S_{21}\left[ S_{32}S_{22}+\frac{\omega_{j2}}{\omega_{j3}}C_{32}C_{22}\right] 
\end{bmatrix}\\
F&=&\begin{bmatrix}
\frac{1}{\omega_{j2}}S_{21}\left[C_{32}C_{22}+\frac{\omega_{j2}}{\omega_{j3}}S_{32}S_{22}\right]+\frac{1}{\omega_{j2}}C_{21}\left[C_{32}S_{22}-\frac{\omega_{j2}}{\omega_{j3}}S_{32}C_{22}\right]\\
\frac{1}{\omega_{j2}}S_{21}\left[S_{32}C_{22}-\frac{\omega_{j2}}{\omega_{j3}}C_{32}S_{22} \right]+\frac{1}{\omega_{j2}}C_{21}\left[S_{32}S_{22}+\frac{\omega_{j2}}{\omega_{j3}}C_{32}C_{22}\right]
\end{bmatrix}.\end{eqnarray}
Simplifying, we can multiply the matrices and obtain
\begin{eqnarray}
{\bf x}_2&=&
\begin{bmatrix}
-C_{21}S_{11}+\frac{\omega_{j1}}{\omega_{j2}}S_{21}C_{11}\\
-S_{21}S_{11}+\frac{\omega_{j1}}{\omega_{j2}}C_{21}C_{11}\\
G+H\\
I+J
\end{bmatrix}
\end{eqnarray}
where
\begin{eqnarray*}
G&=&-S_{11}C_{21}C_{32}C_{22}-\frac{\omega_{j2}}{\omega_{j3}}S_{11}C_{21}S_{32}S_{22}-S_{11}S_{21}C_{32}S_{22}+\frac{\omega_{j2}}{\omega_{j3}}S_{11}S_{21}S_{32}C_{22},\\
H&=&\frac{\omega_{j1}}{\omega_{j2}}C_{11}S_{21}C_{32}C_{22}+\frac{\omega_{j1}}{\omega_{j3}}C_{11}S_{21}S_{32}S_{22}+\frac{\omega_{j1}}{\omega_{j2}}C_{11}C_{21}C_{32}S_{22}-\frac{\omega_{j1}}{\omega_{j3}}C_{11}C_{21}S_{32}C_{22},\\
I&=&-S_{11}C_{21}S_{32}C_{22}+\frac{\omega_{j2}}{\omega_{j3}}S_{11}C_{21}C_{32}S_{22}-S_{11}S_{21}S_{32}S_{22}-\frac{\omega_{j2}}{\omega_{j3}}S_{11}S_{21}C_{32}C_{22},\\
J&=&\frac{\omega_{j1}}{\omega_{j2}}C_{11}S_{21}S_{32}C_{22}-\frac{\omega_{j1}}{\omega_{j3}}C_{11}S_{21}C_{32}S_{22}+\frac{\omega_{j1}}{\omega_{j2}}C_{11}C_{21}S_{32}S_{22}+\frac{\omega_{j1}}{\omega_{j3}}C_{11}C_{21}C_{32}C_{22}.
\end{eqnarray*}

After the null space is found, the values of $\omega_{j1}$ and $\omega_{j2}$ can be found using the $\omega$ relations previously described:
\begin{eqnarray}
\omega_{j1} & = & \sqrt{\omega_{j3}^2 - \alpha_3^2 +\alpha_1^2},\\
\omega_{j2} & = & \sqrt{\omega_{j3}^2 - \alpha_3^2 + \alpha_2^2}.
\end{eqnarray}
Taking all of this information, the eigenfunction
$$V_j(z)=
\left \{
  \begin{array}{ll}
 \sin(\omega_{j1}z)	&	0 \leq z<\pi\rho_1, \\
 A_{j2}\cos(\omega_{j2}z)+B_{j2}\sin(\omega_{j2}z),  &	\pi\rho_1 \leq z<\pi\rho_2\\
  A_{j3}\cos(\omega_{j3}z)+B_{j3}\sin(\omega_{j3}z),	&	\pi\rho_2 \leq z<\pi
    \end{array}
\right.$$
can be constructed for each value of $j$. 


\section{Solving the PDE}
Once the approximate eigenfunctions are constructed, the IBVP
$$\frac{du}{dr} = \pm  i\left(\sqrt{\frac{\partial ^2}{\partial z^2} + \alpha^2I}\right)u, \quad 0<z<\pi , \quad u(z,0)=f(z), \quad u(0,r)=0(\pi,r)=0$$
can be solved through a truncated eigenfunction expansion. Our computed solution is 
$$u(z,r)=\sum_{j=1}^{m}{V_j(z)e^{i\sqrt{\lambda_j}r}}f_j$$
where $f_j=\frac{\langle{V_j,f}\rangle}{\langle{V_j,V_j}\rangle}$ and $\lambda_j=-\omega_{j3}^2+\alpha_3^2$.  The number of terms in the expansion, $m$, can determined through examination of
the coefficients of $V_j(z)$, and taking into account the fact that $e^{i\sqrt{\lambda_j}r}\rightarrow 0$
as $j\rightarrow\infty$.

We will focus on transforming the inner products in $f_j$ to discrete sine transforms. First, the denominator can be handled analytically.
\begin{eqnarray}
\langle{V_j,V_j}\rangle&=&\int_0^\pi{V_j(z)^2}\,dz\nonumber \\
&=&\int_0^{\pi\rho_1}{V_j(z)^2}\,dz+\int_{\pi\rho_1}^{\pi\rho_2}{V_j(z)^2}\,dz+\int_{\pi\rho_2}^{\pi}{V_j(z)^2}
\,dz\nonumber \\
&=&-\frac{1}{4\omega_{j1}}B_{j1}^2(\sin(2\pi\rho_1 \omega_{j1})-2\pi\rho_1\omega_{j1})\nonumber\\
&&-\frac{1}{4\omega_{j2}}[2\pi\omega_{j2} (A_{j2}^2+B_{j2}^2)(\rho_1-\rho_2)+(A_{j2}-B_{j2})(A_{j2}+B_{j2})(\sin(2\pi\rho_1\omega_{j2})-\sin(2\pi\rho_2\omega_{j2}))\nonumber\\&&-2A_{j2}B_{j2}\cos(2\pi\rho_1\omega_{j2})+2A_{j2}B_{j2}\cos(2\pi\rho_2\omega_{j2})]\nonumber\\
&&+\frac{1}{4\omega_{j3}}[-2\pi\omega_{j3}(\rho_2 -1)(A_{j3}^2+B_{j3}^2)+(A_{j3}-B_{j3})(A_{j3}+B_{j3})(\sin(2\pi\omega_{j3})-\sin(2\pi\rho_2\omega_{j3}))\nonumber\\
&&+2A_{j3}B_{j3}\cos(2\pi\rho_2\omega_{j3})-2*A_{j3}B_{j3}\cos(2\pi\omega_{j3})].
\end{eqnarray}
The numerator of $f_j$, $\langle{V_j,f}\rangle$, refers to an inner product and $f(z)$ has the following Fourier sine series:
\begin{eqnarray}
f(z)=\sum_{k=1}^{\infty}{b_k\sin(kz)}.
\end{eqnarray}
Therefore,
\begin{eqnarray}
\langle{V_j,f}\rangle\approx \sum_{k=1}^{m}{b_k\langle{V_j,\sin(k\cdot)}\rangle}, \quad j=1,\ldots,m,
\end{eqnarray}
where $m$ is the number of terms used in our truncated eigenfunction expansion.  We then evaluate $\langle{V_j(z),\sin(k\cdot)}\rangle$:
\begin{eqnarray}
\langle{V_j,\sin(k\cdot)}\rangle&=&\int_{0}^{\pi}{V_j(z)\sin(kz)\,dz}\nonumber\\
&=&\frac{B_{j1}}{k^2-\omega_{j1}^2}[(\omega_{j1}\sin(\pi \rho_1 k)\cos(\pi \rho_1\omega_{j1})-k\cos(\pi \rho_1 k)\sin(\pi \rho_1\omega_{j1})]+\nonumber\\
&&\frac{A_{j2}}{k^2-\omega_{j2}^2}[(\omega_{j2}\sin(\pi k \rho_1)\sin(\pi \omega_{j2} \rho_1)+k\cos(\pi k\rho_1)\cos(\pi\rho_1\omega_{j2})\nonumber\\
&&-\omega_{j2}\sin(\pi k \rho_2)\sin(\pi \rho_2 \omega_{j2})-k\cos(\pi k \rho_2)\cos(\pi\rho_2 \omega_{j2}))]+\nonumber\\
&&\frac{B_{j2}}{k^2-\omega_{j2}^2}[(-\omega_{j2}\sin(\pi k \rho_1)\cos(\pi \omega_{j2} \rho_1)+k\cos(\pi k\rho_1)\sin(\pi\rho_1\omega_{j2})\nonumber\\
&&+\omega_{j2}\sin(\pi k \rho_2)\cos(\pi \rho_2 \omega_{j2})-k\cos(\pi k \rho_2)\sin(\pi\rho_2 \omega_{j2}))]+\nonumber\\
&&\frac{A_{j3}}{k^2-\omega_{j3}^2}[(\omega_{j3}\sin(\pi k \rho_2)\sin(\pi \rho_2 \omega_{j3})+k\cos(\pi k \rho_2)\cos(\pi \rho_2 \omega_{j3})\nonumber\\
&&-\omega_{j3}\sin(\pi k)\sin(\pi\omega_{j3}) -k\cos(\pi k)\cos(\pi\omega_{j3})]+\nonumber\\
&&\frac{B_{j3}}{k^2-\omega_{j3}^2}[(-\omega_{j3}\sin(\pi k \rho_2)\cos(\pi \rho_2 \omega_{j3})+k\cos(\pi k \rho_2)\sin(\pi \rho_2 \omega_{j3})\nonumber\\
&&+\omega_{j3}\sin(\pi k)\cos(\pi\omega_{j3}) -k\cos(\pi k)\sin(\pi\omega_{j3})].
\end{eqnarray}
We then have
$$
{\bf v}=G{\bf b}
$$
where $v_j = \langle V_j, f \rangle$ and $G_{jk} = \langle V_j, \sin(k\cdot) \rangle.$. 

\subsection{Decay Rate of Matrix Entries}
Next, we will examine the asymptotic behavior of the entries of the change-of-basis matrix $G$. For each choice of coefficient $\alpha(z)$, the entries of $G$ behave somewhat differently; however, there are some commonalities that will help us make some inferences about how the structure of the matrix may behave as the number of grid points is increased. 
The coefficients used in this experiment are defined by $\alpha^1=[2,1,2]$ and $\rho^1=[1/3, 2/3]$, $\alpha^2=[1, 3, 5]$ and $\rho^2=[1/4, 3/4]$.  For $k=1,2$, we have
$\alpha(z) = \alpha_i^k$
for $\pi\rho_{i-1}^k \leq z < \pi \rho_i^k$, where for convenience we have defined $\rho_0=0$.
We use $N=127$ grid points in our spatial discretization, and compute the matrices $G_k$, for
$k=1,2$, that are induced by the coefficient values in $\alpha^k$ and interface indicators in $\rho^k$.

\begin{figure}[ht]
\begin{center}
\includegraphics[width=3in]{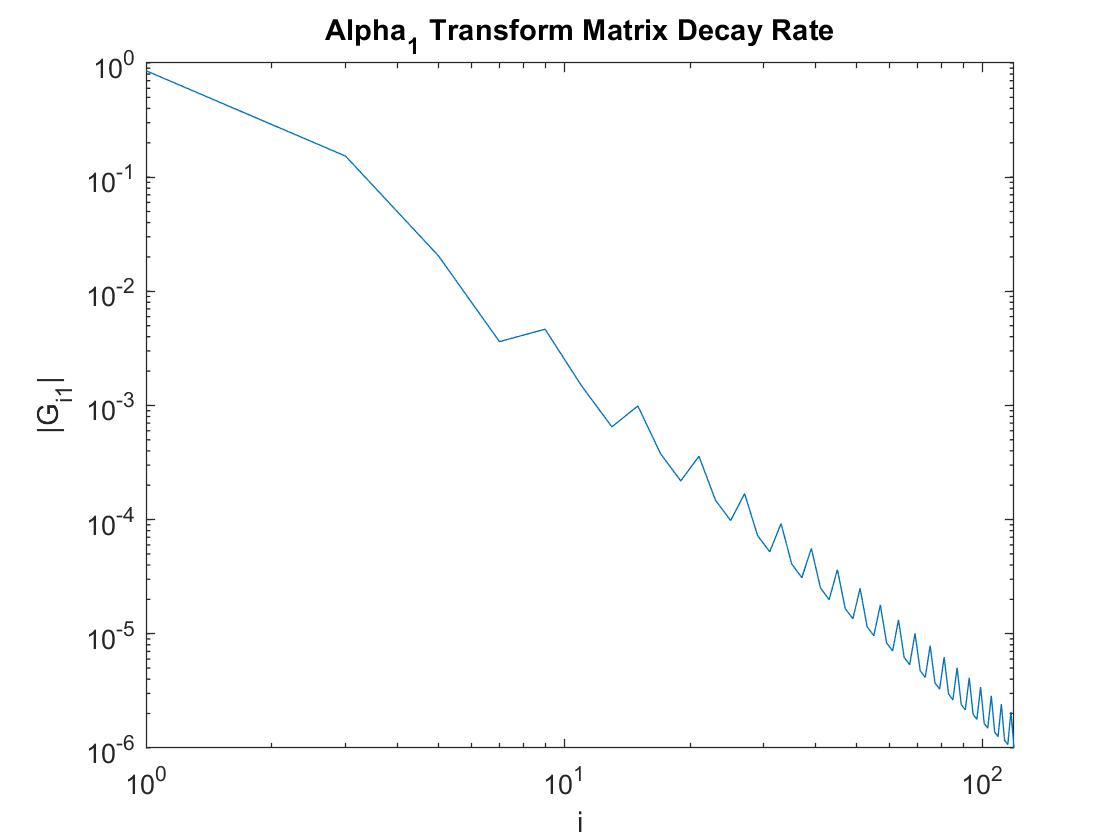} 
\label{alpha1}
\caption{Rate of decay of entries in the first column of the change-of-basis matrix $G_1$.}
\end{center}
\end{figure}

\begin{figure}[ht]
\begin{center}
\includegraphics[width=3in]{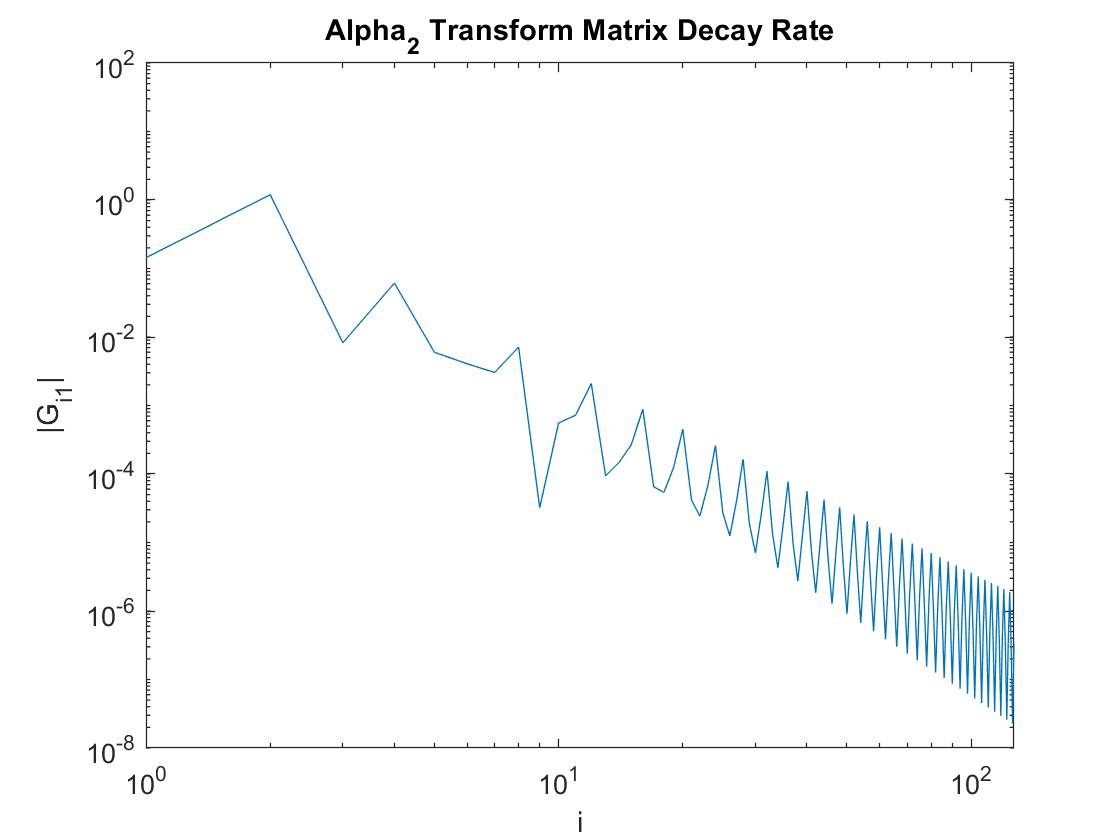} 
\label{alpha2}
\caption{Rate of decay of entries in the first column of the change-of-basis matrix $G_2$.}
\end{center}
\end{figure}

In the case of $G_1$, we observe that alternating entries of each column are zero. Therefore, when studying the rate of decay, we will only consider the nonzero entries. It appears that every sixth entry contributes to a subsequence that is decreasing in magnitude. A plot of the entries from the first column is shown in Figure \ref{alpha1}. The $x$-axis represents the row index of the entry, and the $y$-axis represents the absolute value of the matrix entry. Based on the graph, the order of decay is $j^{-3}$.
Unlike $G_1$, the change-of-basis matrix $G_2$ induced by $\alpha^2$ and $\rho^2$ does not have alternating zero entries.  We observe the pattern that every four entries contributes to a subsequence that is decreasing in magnitude. The plot is shown in Figure \ref{alpha2}. As with $G_1$, the order of decay is $j^{-3}$.
%
%
%

In our algorithm, we computed the coefficients $\langle V_j, f \rangle$ as dot products of grid functions. 
Equivalently, the vector of coefficients is the result of a matrix-vector product of a dense matrix with the grid function for $f(x)$, which requires order $O(N^2)$ operations. It is worthwhile to investigate whether it is more efficient to compute these coefficients using the change-of-basis matrix $G$ in conjunction with a discrete sine transform of $f(x)$, which requires order $O(N\log N)$ operations.  This investigation will be carried out in future work.

\section{Numerical results}


In this section we demonstrate the effectiveness of our approach for solving the IBVP 
(\ref{eq:thePDE}), (\ref{eq:theL}), (\ref{eq:theIC}), with homogeneous Dirichlet boundary conditions.
We compare our approach, labeled ``EigExp" in tables, to the following methods, all of which
use a finite difference discretization $L_N$ of the operator $L$ on a uniform $N$-point grid:
\begin{itemize}
\item Computing the matrix exponential of $iL_N^{1/2}r$ using its Schur Decomposition (labeled ``Schur" in tables),
\item Approximating a matrix-vector product involving $e^{iL_N^{1/2}r}$ using rational Krylov projection, as in \cite{moret}
(labeled ``Krylov" in tables), and
\item The {\sc Matlab} stiff ODE solver {\tt ode15s} (labeled ``{\tt ode15s}" in tables).
\end{itemize}
To measure accuracy, solutions are compared to a reference solution computed using {\tt ode15s} with a
larger number of grid points, $N_{\rm ref} = 8,191$, and the smallest allowed time step.
We calculate the error ${\cal E}_N$ in each computed solution $\tilde{\bf u}^N$ with respect to the reference solution ${\bf u}^{N_{\rm ref}}$ by taking the $\ell_\infty$ norm
of the difference between the computed solution 
and the restriction of the reference solution to the $N$-point grid; that is,
$${\cal E}_N = \frac{\displaystyle{\max_{1\leq i\leq N} \left|u_{m(i)}^{N_{\rm ref}} - \tilde{u}^N_i \right|}}{\displaystyle{\max_{1\leq i\leq N} \left|u_{m(i)}^{N_{\rm ref}} \right|}},$$
where $m(i) = \left( \frac{N_{\rm ref}+1}{N+1}\right) i$ is used to map the values of $\tilde{\bf u}^N$ on its $N$-point grid to the $N_{\rm ref}$-point grid on which the reference solution is computed.

The following numerical results were obtained by using the eigenfunction expansion method described in the previous chapter, which we will call EigExp, to test its efficiency by running it with three different final ranges ($r_f$): 0.01, 0.1, and 1. For each final range, the method was tested using the grid sizes $N=127$, 255, and 511. For each test, the information collected included the time needed to compute the solution, as well as how many terms were needed in the eigenfunction expansion. 

Each numerical experiment was run using either $\sin(2z)$ or a characteristic function as initial data, to test performance with smooth or non-smooth solutions. The piecewise constant coefficient $\alpha(z)$ and the relative locations of the interface within the spatial domain, $\rho_1$ and $\rho_2$, varied between trials as well.  The coefficients used were defined as follows:
\begin{enumerate}
\item $\alpha(z)=\alpha^{1}(z)$, with $\alpha_1^1 =2$, $\alpha_2^1=1$, $\alpha_3^1=2$, 
$\rho_1=\frac{1}{3}$, and $\rho_2=\frac{2}{3}$.
\item $\alpha(z)=\alpha^{2}(z)$, with $\alpha_1^2 =1$, $\alpha_2^2=2$, $\alpha_3^2=3$, 
$\rho_1=\frac{1}{4}$, and $\rho_2=\frac{3}{4}$.
\item $\alpha(z)=\alpha^{3}(z)$, with $\alpha_1^3 =2$, $\alpha_2^3=1$, $\alpha_3^3=3$, 
$\rho_1=\frac{1}{4}$, and $\rho_2=\frac{1}{2}$.
\end{enumerate}

The first six tables below report the time required for each method to compute the solution for
each problem.  For EigExp, ``terms'' refers to the number of terms, $m$, included in the truncated
eigenfunction expansion of the solution.  For Krylov, ``terms'' refers to the number of Krylov projection
steps taken to achieve convergence.  For {\tt ode15}, $\Delta t_{\rm avg}$ is the average time step
taken. From these tables, we can make the following observations:
\begin{itemize}
\item In terms of efficiency, only rational Krylov projection was competitive with our eigenfunction
expansion algorithm.
\item For smooth ``initial'' data, Krylov projection was more efficient for smaller $r_f$, but this edge disappeared as $r_f$ increased, with the eigenfunction expansion requiring substantially fewer terms.
\item For discontinuous initial data, eigenfunction expansion was more efficient regardless of $r_f$.
\item Eigenfunction expansion was the most accurate of all of the methods, except
in the case of discontinuous initial data and medium-to-large $r_f$, in which case either Schur
Decomposition or Krylov projection was slightly more accurate.
\end{itemize}

\begin{center}
\begin{tabular}{|l|c|c|c|c|c|c|c|c|} \hline
\multicolumn{9}{|c|}{$\alpha^1$, $\sin(2z)$ } \\ \hline
\multicolumn{2}{|c|}{} & \multicolumn{2}{|c|}{EigExp} & Schur & \multicolumn{2}{|c|}{Krylov} & \multicolumn{2}{|c|}{{\tt ode15s}} \\ \hline
$r_f$ & $N$ & time & terms & time & time & terms & time & $\Delta t_{\mathrm{avg}}$ \\ \hline
0.01 & 127 & 0.151 & 120 & 0.123 & 0.099 & 16 & 0.102 & 1.000e-03\\
 & 255 & 0.055 & 114 & 0.295 & 0.019 & 11 & 0.069 & 1.000e-03\\
 & 511 & 0.027 & 114 & 2.020 & 0.007 & 11 & 0.152 & 1.000e-03\\ \hline
0.1 & 127 & 0.015 & 44 & 0.105 & 0.042 & 10 & 0.173 & 1.000e-02\\
 & 255 & 0.012 & 44 & 0.338 & 0.022 & 9 & 0.046 & 1.000e-02\\
 & 511 & 0.011 & 44 & 1.943 & 0.006 & 10 & 0.191 & 1.000e-02\\ \hline
1 & 127 & 0.007 & 14 & 0.103 & 0.011 & 6 & 0.014 & 6.667e-02\\
 & 255 & 0.006 & 14 & 0.255 & 0.004 & 6 & 0.056 & 6.667e-02\\
 & 511 & 0.006 & 14 & 2.235 & 0.009 & 6 & 0.323 & 6.667e-02\\ \hline
\end{tabular}

\begin{tabular}{|l|c|c|c|c|c|c|c|c|} \hline
\multicolumn{9}{|c|}{$\alpha^1$, characteristic function}\\ \hline
\multicolumn{2}{|c|}{} & \multicolumn{2}{|c|}{EigExp} & Schur & \multicolumn{2}{|c|}{Krylov} & \multicolumn{2}{|c|}{{\tt ode15s}} \\ \hline
$r_f$ & $N$ & time & terms & time & time & terms & time & $\Delta t_{\mathrm{avg}}$ \\ \hline
0.01 & 127 & 0.015 & 128 & 0.098 & 0.095 & 33 & 0.021 & 5.882e-04\\
 & 255 & 0.028 & 256 & 0.269 & 0.069 & 41 & 0.059 & 5.263e-04\\
 & 511 & 0.131 & 512 & 1.836 & 0.207 & 58 & 0.249 & 4.762e-04\\ \hline
0.1 & 127 & 0.023 & 104 & 0.100 & 0.028 & 27 & 0.016 & 3.571e-03\\
 & 255 & 0.024 & 97 & 0.284 & 0.024 & 28 & 0.061 & 2.941e-03\\
 & 511 & 0.025 & 97 & 1.799 & 0.051 & 28 & 0.458 & 2.500e-03\\ \hline
1 & 127 & 0.007 & 17 & 0.073 & 0.004 & 12 & 0.025 & 2.041e-02\\
 & 255 & 0.003 & 17 & 0.271 & 0.005 & 12 & 0.085 & 1.818e-02\\
 & 511 & 0.007 & 17 & 1.724 & 0.018 & 12 & 0.600 & 1.639e-02\\ \hline
\end{tabular}

\begin{tabular}{|l|c|c|c|c|c|c|c|c|} \hline
\multicolumn{9}{|c|}{$\alpha^2$, $\sin(2z)$ }\\ \hline
\multicolumn{2}{|c|}{} & \multicolumn{2}{|c|}{EigExp} & Schur & \multicolumn{2}{|c|}{Krylov} & \multicolumn{2}{|c|}{{\tt ode15s}} \\ \hline
$r_f$ & $N$ & time & terms & time & time & terms & time & $\Delta t_{\mathrm{avg}}$ \\ \hline
0.01 & 127 & 0.051 & 128 & 0.081 & 0.017 & 21 & 0.021 & 1.000e-03\\
 & 255 & 0.013 & 127 & 0.313 & 0.022 & 26 & 0.032 & 1.000e-03\\
 & 511 & 0.047 & 127 & 1.777 & 0.035 & 26 & 0.140 & 1.000e-03\\ \hline
0.1 & 127 & 0.010 & 47 & 0.105 & 0.010 & 17 & 0.012 & 1.000e-02\\
 & 255 & 0.008 & 47 & 0.268 & 0.012 & 17 & 0.032 & 1.000e-02\\
 & 511 & 0.017 & 47 & 1.745 & 0.020 & 17 & 0.174 & 1.000e-02\\ \hline
1 & 127 & 0.006 & 14 & 0.133 & 0.005 & 8 & 0.014 & 5.882e-02\\
 & 255 & 0.003 & 14 & 0.271 & 0.003 & 8 & 0.057 & 5.263e-02\\
 & 511 & 0.006 & 14 & 1.751 & 0.006 & 8 & 0.258 & 5.000e-02\\ \hline
\end{tabular}

\begin{tabular}{|l|c|c|c|c|c|c|c|c|} \hline
\multicolumn{9}{|c|}{$\alpha^2$, characteristic function}\\ \hline
\multicolumn{2}{|c|}{} & \multicolumn{2}{|c|}{EigExp} & Schur & \multicolumn{2}{|c|}{Krylov} & \multicolumn{2}{|c|}{{\tt ode15s}} \\ \hline
$r_f$ & $N$ & time & terms & time & time & terms & time & $\Delta t_{\mathrm{avg}}$ \\ \hline
0.01 & 127 & 0.021 & 128 & 0.132 & 0.031 & 33 & 0.011 & 5.882e-04\\
 & 255 & 0.030 & 256 & 0.259 & 0.053 & 41 & 0.046 & 5.263e-04\\
 & 511 & 0.089 & 512 & 1.805 & 0.141 & 58 & 0.255 & 4.762e-04\\ \hline
0.1 & 127 & 0.019 & 104 & 0.083 & 0.022 & 27 & 0.015 & 3.448e-03\\
 & 255 & 0.016 & 97 & 0.381 & 0.055 & 28 & 0.099 & 2.941e-03\\
 & 511 & 0.022 & 97 & 2.049 & 0.034 & 28 & 0.480 & 2.500e-03\\ \hline
1 & 127 & 0.005 & 17 & 0.098 & 0.006 & 12 & 0.031 & 2.000e-02\\
 & 255 & 0.006 & 17 & 0.341 & 0.008 & 12 & 0.138 & 1.786e-02\\
 & 511 & 0.006 & 17 & 1.867 & 0.010 & 12 & 0.573 & 1.613e-02\\ \hline
\end{tabular}

\begin{tabular}{|l|c|c|c|c|c|c|c|c|} \hline
\multicolumn{9}{|c|}{$\alpha^3$, $\sin(2z)$ }\\ \hline
\multicolumn{2}{|c|}{} & \multicolumn{2}{|c|}{EigExp} & Schur & \multicolumn{2}{|c|}{Krylov} & \multicolumn{2}{|c|}{{\tt ode15s}} \\ \hline
$r_f$ & $N$ & time & terms & time & time & terms & time & $\Delta t_{\mathrm{avg}}$ \\ \hline
0.01 & 127 & 0.036 & 110 & 0.078 & 0.010 & 21 & 0.008 & 1.000e-03\\
 & 255 & 0.025 & 107 & 0.317 & 0.018 & 23 & 0.031 & 1.000e-03\\
 & 511 & 0.024 & 106 & 1.894 & 0.020 & 22 & 0.140 & 1.000e-03\\ \hline
0.1 & 127 & 0.009 & 43 & 0.096 & 0.005 & 13 & 0.008 & 1.000e-02\\
 & 255 & 0.014 & 43 & 0.297 & 0.005 & 13 & 0.030 & 1.000e-02\\
 & 511 & 0.010 & 43 & 1.700 & 0.007 & 13 & 0.173 & 1.000e-02\\ \hline
1 & 127 & 0.004 & 14 & 0.086 & 0.004 & 9 & 0.023 & 4.762e-02\\
 & 255 & 0.006 & 14 & 0.302 & 0.004 & 9 & 0.056 & 4.545e-02\\
 & 511 & 0.011 & 14 & 1.841 & 0.005 & 9 & 0.299 & 4.762e-02\\ \hline
\end{tabular}

\begin{tabular}{|l|c|c|c|c|c|c|c|c|} \hline
\multicolumn{9}{|c|}{$\alpha^3$, characteristic function}\\ \hline
\multicolumn{2}{|c|}{} & \multicolumn{2}{|c|}{EigExp} & Schur & \multicolumn{2}{|c|}{Krylov} & \multicolumn{2}{|c|}{{\tt ode15s}} \\ \hline
$r_f$ & $N$ & time & terms & time & time & terms & time & $\Delta t_{\mathrm{avg}}$ \\ \hline
0.01 & 127 & 0.014 & 128 & 0.080 & 0.028 & 33 & 0.010 & 5.882e-04\\
 & 255 & 0.030 & 256 & 0.289 & 0.050 & 42 & 0.042 & 5.263e-04\\
 & 511 & 0.066 & 512 & 1.764 & 0.161 & 58 & 0.257 & 4.762e-04\\ \hline
0.1 & 127 & 0.016 & 104 & 0.116 & 0.037 & 27 & 0.031 & 3.448e-03\\
 & 255 & 0.012 & 97 & 0.281 & 0.044 & 28 & 0.093 & 2.941e-03\\
 & 511 & 0.021 & 97 & 2.020 & 0.052 & 29 & 0.441 & 2.500e-03\\ \hline
1 & 127 & 0.006 & 17 & 0.084 & 0.005 & 12 & 0.033 & 2.041e-02\\
 & 255 & 0.009 & 17 & 0.352 & 0.011 & 12 & 0.160 & 1.818e-02\\
 & 511 & 0.007 & 17 & 1.812 & 0.008 & 12 & 0.646 & 1.639e-02\\ \hline
\end{tabular}

\begin{tabular}{|l|c|c|c|c|c|c|c|c|} \hline
\multicolumn{9}{|c|}{$\alpha^1$, $\sin(2z)$ }\\ \hline
\multicolumn{2}{|c|}{} & \multicolumn{2}{|c|}{EigExp} & Schur & \multicolumn{2}{|c|}{Krylov} & \multicolumn{2}{|c|}{{\tt ode15s}} \\ \hline
$r_f$ & $N$ & error & terms & error & error & terms & error & $\Delta t_{\mathrm{avg}}$ \\ \hline
0.01 & 127 & 1.4e-04 & 120 & 1.7e-04 & 1.7e-04 & 16 & 1.6e-04 & 1.000e-03\\
 & 255 & 6.5e-06 & 114 & 9.9e-05 & 9.5e-05 & 11 & 1.0e-04 & 1.000e-03\\
 & 511 & 7.4e-06 & 114 & 4.5e-05 & 5.5e-05 & 11 & 4.0e-05 & 1.000e-03\\ \hline
0.1 & 127 & 7.3e-05 & 44 & 1.1e-03 & 1.1e-03 & 10 & 9.6e-04 & 1.000e-02\\
 & 255 & 1.8e-05 & 44 & 5.9e-04 & 6.0e-04 & 9 & 6.9e-04 & 1.000e-02\\
 & 511 & 1.8e-05 & 44 & 2.6e-04 & 2.6e-04 & 10 & 1.5e-04 & 1.000e-02\\ \hline
1 & 127 & 1.1e-04 & 14 & 3.1e-03 & 3.1e-03 & 6 & 3.0e-03 & 6.667e-02\\
 & 255 & 1.1e-04 & 14 & 1.9e-03 & 1.9e-03 & 6 & 2.1e-03 & 6.667e-02\\
 & 511 & 1.1e-04 & 14 & 8.1e-04 & 8.1e-04 & 6 & 6.7e-04 & 6.667e-02\\ \hline
\end{tabular}

\begin{tabular}{|l|c|c|c|c|c|c|c|c|} \hline
\multicolumn{9}{|c|}{$\alpha^1$, characteristic function}\\ \hline
\multicolumn{2}{|c|}{} & \multicolumn{2}{|c|}{EigExp} & Schur & \multicolumn{2}{|c|}{Krylov} & \multicolumn{2}{|c|}{{\tt ode15s}} \\ \hline
$r_f$ & $N$ & error & terms & error & error & terms & error & $\Delta t_{\mathrm{avg}}$ \\ \hline
0.01 & 127 & 2.8e-01 & 128 & 3.0e-01 & 3.0e-01 & 33 & 3.0e-01 & 5.882e-04\\
 & 255 & 1.8e-01 & 256 & 2.0e-01 & 2.0e-01 & 41 & 2.0e-01 & 5.263e-04\\
 & 511 & 9.2e-02 & 512 & 1.0e-01 & 1.0e-01 & 58 & 1.0e-01 & 4.762e-04\\ \hline
0.1 & 127 & 4.3e-02 & 104 & 4.1e-02 & 4.1e-02 & 27 & 8.8e-02 & 3.571e-03\\
 & 255 & 2.1e-02 & 97 & 2.0e-02 & 2.0e-02 & 28 & 8.0e-02 & 2.941e-03\\
 & 511 & 1.0e-02 & 97 & 9.7e-03 & 9.7e-03 & 28 & 7.9e-02 & 2.500e-03\\ \hline
1 & 127 & 2.5e-02 & 17 & 1.1e-02 & 1.1e-02 & 12 & 6.3e-01 & 2.041e-02\\
 & 255 & 1.2e-02 & 17 & 6.2e-03 & 6.2e-03 & 12 & 6.3e-01 & 1.818e-02\\
 & 511 & 6.1e-03 & 17 & 2.7e-03 & 2.7e-03 & 12 & 6.2e-01 & 1.639e-02\\ \hline
\end{tabular}

\begin{tabular}{|l|c|c|c|c|c|c|c|c|} \hline
\multicolumn{9}{|c|}{$\alpha^2$, $\sin(2z)$ }\\ \hline
\multicolumn{2}{|c|}{} & \multicolumn{2}{|c|}{EigExp} & Schur & \multicolumn{2}{|c|}{Krylov} & \multicolumn{2}{|c|}{{\tt ode15s}} \\ \hline
$r_f$ & $N$ & error & terms & error & error & terms & error & $\Delta t_{\mathrm{avg}}$ \\ \hline
0.01 & 127 & 1.2e-05 & 128 & 6.9e-04 & 6.9e-04 & 21 & 3.3e-02 & 1.000e-03\\
 & 255 & 1.4e-05 & 127 & 3.7e-04 & 3.6e-04 & 26 & 3.2e-02 & 1.000e-03\\
 & 511 & 1.6e-05 & 127 & 1.8e-04 & 1.8e-04 & 26 & 3.2e-02 & 1.000e-03\\ \hline
0.1 & 127 & 6.1e-05 & 47 & 4.0e-03 & 4.0e-03 & 17 & 3.2e-01 & 1.000e-02\\
 & 255 & 6.1e-05 & 47 & 2.0e-03 & 2.0e-03 & 17 & 3.2e-01 & 1.000e-02\\
 & 511 & 6.1e-05 & 47 & 9.4e-04 & 9.4e-04 & 17 & 3.2e-01 & 1.000e-02\\ \hline
1 & 127 & 3.2e-04 & 14 & 2.1e-02 & 2.1e-02 & 8 & 2.0e+00 & 5.882e-02\\
 & 255 & 3.2e-04 & 14 & 1.0e-02 & 1.0e-02 & 8 & 2.0e+00 & 5.263e-02\\
 & 511 & 3.2e-04 & 14 & 4.9e-03 & 4.9e-03 & 8 & 2.0e+00 & 5.000e-02\\ \hline
\end{tabular}

\begin{tabular}{|l|c|c|c|c|c|c|c|c|} \hline
\multicolumn{9}{|c|}{$\alpha^2$, characteristic function}\\ \hline
\multicolumn{2}{|c|}{} & \multicolumn{2}{|c|}{EigExp} & Schur & \multicolumn{2}{|c|}{Krylov} & \multicolumn{2}{|c|}{{\tt ode15s}} \\ \hline
$r_f$ & $N$ & error & terms & error & error & terms & error & $\Delta t_{\mathrm{avg}}$ \\ \hline
0.01 & 127 & 2.8e-01 & 128 & 3.0e-01 & 3.0e-01 & 33 & 3.0e-01 & 5.882e-04\\
 & 255 & 1.8e-01 & 256 & 2.0e-01 & 2.0e-01 & 41 & 2.0e-01 & 5.263e-04\\
 & 511 & 9.2e-02 & 512 & 1.0e-01 & 1.0e-01 & 58 & 1.0e-01 & 4.762e-04\\ \hline
0.1 & 127 & 4.6e-02 & 104 & 4.2e-02 & 4.2e-02 & 27 & 1.3e-01 & 3.448e-03\\
 & 255 & 2.2e-02 & 97 & 2.0e-02 & 2.0e-02 & 28 & 1.2e-01 & 2.941e-03\\
 & 511 & 1.1e-02 & 97 & 9.8e-03 & 9.8e-03 & 28 & 1.2e-01 & 2.500e-03\\ \hline
1 & 127 & 2.8e-02 & 17 & 1.5e-02 & 1.5e-02 & 12 & 9.1e-01 & 2.000e-02\\
 & 255 & 1.4e-02 & 17 & 7.2e-03 & 7.2e-03 & 12 & 9.0e-01 & 1.786e-02\\
 & 511 & 6.7e-03 & 17 & 3.5e-03 & 3.5e-03 & 12 & 9.0e-01 & 1.613e-02\\ \hline
\end{tabular}

\begin{tabular}{|l|c|c|c|c|c|c|c|c|} \hline
\multicolumn{9}{|c|}{$\alpha^3$, $\sin(2z)$ }\\ \hline
\multicolumn{2}{|c|}{} & \multicolumn{2}{|c|}{EigExp} & Schur & \multicolumn{2}{|c|}{Krylov} & \multicolumn{2}{|c|}{{\tt ode15s}} \\ \hline
$r_f$ & $N$ & error & terms & error & error & terms & error & $\Delta t_{\mathrm{avg}}$ \\ \hline
0.01 & 127 & 1.1e-05 & 110 & 6.6e-04 & 6.6e-04 & 21 & 3.9e-02 & 1.000e-03\\
 & 255 & 1.6e-05 & 107 & 3.4e-04 & 3.4e-04 & 23 & 3.9e-02 & 1.000e-03\\
 & 511 & 1.7e-05 & 106 & 1.7e-04 & 1.7e-04 & 22 & 3.9e-02 & 1.000e-03\\ \hline
0.1 & 127 & 7.0e-05 & 43 & 4.5e-03 & 4.5e-03 & 13 & 3.9e-01 & 1.000e-02\\
 & 255 & 7.0e-05 & 43 & 2.2e-03 & 2.2e-03 & 13 & 3.9e-01 & 1.000e-02\\
 & 511 & 7.0e-05 & 43 & 1.1e-03 & 1.1e-03 & 13 & 3.9e-01 & 1.000e-02\\ \hline
1 & 127 & 4.1e-04 & 14 & 2.0e-02 & 2.0e-02 & 9 & 1.4e+00 & 4.762e-02\\
 & 255 & 4.1e-04 & 14 & 1.0e-02 & 1.0e-02 & 9 & 1.4e+00 & 4.545e-02\\
 & 511 & 4.1e-04 & 14 & 5.0e-03 & 5.1e-03 & 9 & 1.4e+00 & 4.762e-02\\ \hline
\end{tabular}

\begin{tabular}{|l|c|c|c|c|c|c|c|c|} \hline
\multicolumn{9}{|c|}{$\alpha^3$, characteristic function}\\ \hline
\multicolumn{2}{|c|}{} & \multicolumn{2}{|c|}{EigExp} & Schur & \multicolumn{2}{|c|}{Krylov} & \multicolumn{2}{|c|}{{\tt ode15s}} \\ \hline
$r_f$ & $N$ & error & terms & error & error & terms & error & $\Delta t_{\mathrm{avg}}$ \\ \hline
0.01 & 127 & 2.8e-01 & 128 & 3.0e-01 & 3.0e-01 & 33 & 3.0e-01 & 5.882e-04\\
 & 255 & 1.8e-01 & 256 & 2.0e-01 & 2.0e-01 & 42 & 2.0e-01 & 5.263e-04\\
 & 511 & 9.2e-02 & 512 & 1.0e-01 & 1.0e-01 & 58 & 1.0e-01 & 4.762e-04\\ \hline
0.1 & 127 & 4.4e-02 & 104 & 4.3e-02 & 4.3e-02 & 27 & 7.6e-02 & 3.448e-03\\
 & 255 & 2.2e-02 & 97 & 2.1e-02 & 2.1e-02 & 28 & 7.6e-02 & 2.941e-03\\
 & 511 & 1.1e-02 & 97 & 1.0e-02 & 1.0e-02 & 29 & 7.6e-02 & 2.500e-03\\ \hline
1 & 127 & 2.6e-02 & 17 & 2.9e-02 & 2.9e-02 & 12 & 5.2e-01 & 2.041e-02\\
 & 255 & 1.2e-02 & 17 & 1.4e-02 & 1.4e-02 & 12 & 5.2e-01 & 1.818e-02\\
 & 511 & 5.7e-03 & 17 & 7.0e-03 & 7.0e-03 & 12 & 5.2e-01 & 1.639e-02\\ \hline
\end{tabular}
\end{center}

\section{Generalization to Arbitrary Piecewise Constant Coefficients}
To this point, we have focused on the three-piece case.  In this section, we will show how 
our approach to computing the eigenvalues and eigenfunctions can be generalized to the case
of a piecewise constant coefficient $\alpha(z)$ consisting of $n$ pieces.  The boundary 
and continuity conditions lead to a homogeneous system of equations for the coefficients
$A_{jk}$ and $B_{jk}$, $k=1,\ldots,n$, for the $j$th eigenfunction.  The coefficient matrix of this system is

$$M_j=\left[
\begin{array}{cccccc}
0 & 0 & 0 &  \cdots & \cdots &  {\bf b}_{1n}^T \\
{\bf c}_{11} & E_{12} & 0 & \cdots & \cdots& 0\\
0 & E_{22} & E_{23} & 0 & \cdots & 0\\
0 & 0 & E_{33} & \ddots & 0 & 0\\
\vdots & & & \ddots & \ddots &0\\
&&&&E_{n-1,n-1} & E_{n-1,n}\\
\end{array}
\right]$$
where 
$$E_{ij}=\left[
\begin{array}{cc}
-C_{ij} & -S_{ij}\\
\omega_jS_{ij} & -\omega_jC_{ij}\\
\end{array}
\right]$$
and, as before, $C_{ik}=\cos(\omega_{ji}\rho_k\pi)$ and $S_{ik}=\sin(\omega_{ji}\rho_k\pi)$.
The vector
$${\bf b}_{1n}^T=\left[ \begin{array}{cc} C_{1n} & S_{1n}\end{array}\right]$$ comes from the boundary conditions,
while $${\bf c}_{11}=\left[ \begin{array}{c} S_{11} \\ \omega_{j1}S_{11}\end{array}\right]$$ comes from the continuity conditions.

By analogy with the three-piece case, we require an algorithm to find a value of $\omega_{jn}$ for which $M_j$ is singular. We have
$$\det(M_j)=\det(-B_j D_j^{-1}C_j)\det(D_j)$$
where the blocks $B_j$, $C_j$, $D_j$ are given by
$$B_j=\left[
\begin{array}{cccc}
0 &   \cdots & 0 & {\bf b}_{1n}^T\\
\end{array} \right], \quad 
C_j=\left[
\begin{array}{c}
{\bf c}_{11} \\
0\\
\vdots\\
0
\end{array} \right],
\quad
D_j=\left[
\begin{array}{ccccc}
E_{12} & 0 & \cdots & \cdots& 0\\
E_{22} & E_{23} & 0 & \cdots & 0\\
0 & E_{33} & \ddots & 0 & 0\\
0 & 0 & \ddots & \ddots &0\\
0&0&0&E_{n-1,n-1} & E_{n-1,n}\\
\end{array} \right].$$
As in the three-piece case, the block lower triangular structure of $D_j$ can be used to easily
prove that $D_j$ is nonsingular.

To obtain initial guesses for eigenvalues for the indices $j$ such that $\sin(jz)$ is not an accurate approximation of an eigenfunction, we again construct the matrix $L_0$ with entries
$[L_0]_{ij} =\langle \phi_i, L\phi_j \rangle$. 
For the $i\neq j$ case, we have
\begin{eqnarray}
\left[ L_0 \right]_{ij}&=&\displaystyle{-\left(\frac{j^2}{\pi}\right)\frac{(i+j)\sin(i\pi-j\pi)-(i-j)\sin(i\pi+j\pi)}{i^2-j^2}}\nonumber\\
&+&\frac{\alpha_1^2}{\pi}\left[ \frac{(i+j)\sin(\rho_1(i\pi-j\pi))-(i-j)\sin(\rho_1(i\pi+j\pi))}{i^2-j^2}\right]\nonumber\\
&+&\sum_{v=2}^{n-1}\left\{\frac{\alpha_v^2}{\pi}\left[ \frac{(i+j)\sin(\rho_{v}(i\pi-j\pi))-(i+j)\sin(\rho_{v-1}(i\pi-j\pi))}{i^2-j^2}\right.\right. \nonumber\\
&+&\left.\left.\frac{-(i-j)\sin(\rho_v(i\pi+j\pi))+(i-j)\sin(\rho_{v-1}(i\pi+j\pi))}{i^2-j^2}\right]\right\}\nonumber\\
&+&\frac{\alpha_n^2}{\pi}\left[ \frac{-(i+j)\sin(\rho_{n-1}(i\pi-j\pi))+(i-j)\sin(\rho_{n-1}(i\pi+j\pi))}{i^2-j^2}\right],
\end{eqnarray}
while for the $i=j$ case, 
\begin{eqnarray}
\left[L_0 \right]_{ij}&=&\frac{\langle{ \phi_{j},L\phi_{j} }\rangle}{\langle{\phi_{j},\phi_{j} }\rangle}=-j^2+\frac{2}{\pi}\int_{0}^{\pi}\alpha(z)^2\sin^2\left(jz\right)\,dz\nonumber \\
&=&-j^2+\frac{2}{\pi}B 
\end{eqnarray}
where
\begin{eqnarray*}
B&=& \alpha_{1}^2\left[\frac{-\sin(\pi\rho_{1}j)+\pi\rho_{1}j}{4j}\right] +\nonumber\\
&&\sum_{v=2}^{n-1}\alpha_{v}^2\left[\frac{-\sin(\pi\rho_{v}j)+\pi\rho_{v}j+\sin(\pi\rho_{v-1}j)-\pi\rho_{v-1}j}{4j}\right]  + \alpha_{n}^2\left[\pi + \frac{\sin(\pi\rho_{n-1}j)-\pi\rho_{n-1}j}{4j}\right].
\end{eqnarray*}

As before, to determine the indices $j$ for which $\sin(jz)$ is an accurate approximate eigenfunction,
we use the error formula
$$E_{j}=\frac{\langle{L\phi_{j},L\phi_{j}}\rangle -\mu_{j}^2 \langle{ \phi_{j},\phi_{j} }\rangle}{\langle{L\phi_{j},L\phi_{j} }\rangle}$$
where $$\mu_{j}=\frac{\langle{ \phi_{j},L\phi_{j} }\rangle}{\langle{\phi_{j},\phi_{j} }\rangle}.$$
We obtain
\begin{eqnarray}
E_{j}&=&
\frac{\displaystyle{ A_j-\frac{2}{\pi}B_j^2}}{\displaystyle{\left(\frac{j^4\pi}{2}\right)-2j^2B_j+A_j}}
\end{eqnarray}
where
\begin{eqnarray}
A_j &=& \alpha_{1}^4 \int_{0}^{\pi\rho_{1}} \sin^2(jz)\, dz + \sum_{v=2}^{n-1}\left\{\alpha_{v}^4 \int_{\pi\rho_{v-1}}^{\pi\rho_{v}} \sin^2(jz)\, dz\right\} + \alpha_{n}^4 \int_{\pi\rho_{n-1}}^{\pi} \sin^2(jz)\,dz \\
B_j &=& \alpha_{1}^2 \int_{0}^{\pi\rho_{1}} \sin^2(jz)\, dz + \sum_{v=2}^{n-1}\left\{\alpha_{v}^2 \int_{\pi\rho_{v-1}}^{\pi\rho_{v}} \sin^2(jz)\, dz\right\} + \alpha_{n}^2 \int_{\pi\rho_{n-1}}^{\pi} \sin^2(jz)\, dz
\end{eqnarray}
Evaluating these integrals yields
\begin{eqnarray}
A_j&=&\alpha_{1}^4\left[-\frac{\sin(2\pi\rho_{1}j)-2\pi\rho_{1}j}{4j}\right]\nonumber\\
&& +\sum_{v=2}^{n-1}\left\{\alpha_{v}^4\left[-\frac{\sin(2\pi\rho_{v}j)-2\pi\rho_{v}j-\sin(2\pi\rho_{v-1}j)+2\pi\rho_{v-1}j}{4j}\right]\right\}\nonumber \\
&& + \alpha_{n}^4\left[\frac{\pi}{2} + \frac{\sin(2\pi\rho_{n-1}j)-2\pi\rho_{n-1}j}{4j}\right],  \\
B_j&=& \alpha_{1}^2\left[-\frac{\sin(2\pi\rho_{1}j)-2\pi\rho_{1}j}{4j}\right]\nonumber\\
&& +\sum_{v=2}^{n-1}\left\{\alpha_{v}^2\left[-\frac{\sin(2\pi\rho_{v}j)-2\pi\rho_{v}j-\sin(2\pi\rho_{v-1}j)+2\pi\rho_{v-1}j}{4j}\right]\right\}\nonumber \\
&& + \alpha_{n}^2\left[\frac{\pi}{2} + \frac{\sin(2\pi\rho_{n-1}j)-2\pi\rho_{n-1}j}{4j}\right]. 
\end{eqnarray}

Once the initial guesses for $\omega_{jn}$ are obtained and passed to the Secant Method so that accurate approximations of $\omega_{jn}$ values can be computed, then the matrix $M_j$ can be completed using this $\omega_{jn}$ value.  The null space of $M_j$ is given by the span of
$${\bf x}_j = \left[ \begin{array}{c} 1 \\ -D_j^{-1} C_j \end{array} \right],$$
which can be computed efficiently due to the block lower bidiagonal structure of $D_j$, in which
each block is $2\times 2$.
The remaining values of $\omega_{ji}$, for $i=1,\ldots,n-1$, are found using the relation
\begin{eqnarray}
\omega_{ji} & = & \sqrt{\omega_{jn}^2 - \alpha_n^2 +\alpha_i^2}.
\end{eqnarray}
Finally, the $j$th eigenfunction can be constructed as follows:
$$V_j(z)=
\left \{
  \begin{array}{ll}
  B_{j1}\sin(\omega_{j1}z)	&	0 \leq x<\pi\rho_1, \\
 A_{j2}\cos(\omega_{j2}z)+B_{j2}\sin(\omega_{j2}z)  &	\pi\rho_1 \leq x<\pi\rho_2,\\
\vdots\\
  A_{jn}\cos(\omega_{jn}z)+B_{jn}\sin(\omega_{jn}z)	&	\pi\rho_{n-1} \leq x<\pi,
    \end{array}
\right.$$
where $B_{j1}, A_{j2}, B_{j2}, \ldots, A_{jn}, B_{jn}$ are the elements of the null space vector ${\bf x}_j$.

\section{Conclusion}
The goal of this project was to develop an algorithm to efficiently solve an extra-wide angle parabolic equation with a piecewise constant coefficient
by computing a truncated eigenfunction expansion of the solution. 
We found that the function $\sin(jz)$ was an accurate approximation of the eigenfunctions associated with larger eigenvalues; however, those associated with smaller eigenvalues had to be approximated by a linear combination of low-frequency sines.  These initial approximations to eigenfunctions yielded
initial guesses for the Secant Method. The Secant Method was used to find the exact eigenvalues for each piece, and then the eigenfunctions can be constructed from these. 
Numerical experiments demonstrated the accuracy and efficiency of our approach, relative to existing
numerical methods.

Ongoing work involves the approximate diagonalization of the operator $L$ using the Fokas
transform \cite{fokas}--in particular, exploring the relationship between the 
forward Fokas transform and inner products with eigenfunctions.  An related direction for future
work entails generalization of our algorithm to smoothly varying $\alpha^2$ by taking the limit as
the number of pieces approaches infinity, as was accomplished through the Universal Transform
Method \cite{farkas}.  For problems in which the sound speed is range-dependent, eigenfunctions
computed using our algorithm, using the average value of the sound speed with respect to range,
can be used as basis functions for a Krylov subspace spectral method \cite{nmpde}, which is most
effective when the basis functions are approximate eigenfunctions.






\bibliographystyle{elsarticle-num} 


\bibliography{mybibfile}

%
%
%
\end{document}